\documentclass[10pt,11pt]{article}%
\usepackage[dvips]{graphicx}
\usepackage{amsmath}
\usepackage{amsfonts}
\usepackage{amssymb}%
\providecommand{\U}[1]{\protect\rule{.1in}{.1in}}
\providecommand{\U}[1]{\protect\rule{.1in}{.1in}}
\providecommand{\U}[1]{\protect\rule{.1in}{.1in}}
\providecommand{\U}[1]{\protect\rule{.1in}{.1in}}
\providecommand{\U}[1]{\protect\rule{.1in}{.1in}}
\newtheorem{theorem}{Theorem}[subsection]

\newtheorem{proposition}[theorem]{Proposition}

\newtheorem{fact}[theorem]{Fact}
\newtheorem{property}[theorem]{Property}

\def \E{\mathbb{E}}

\begin{document}

\title{Return Probabilities for the Reflected Random Walk on $\mathbb N_0$}
\date{22/5/2012}
\author{Rim Essifi \& Marc Peign\'e}
\maketitle

\begin{abstract} Let $(Y_n)$ be a sequence of i.i.d. $\mathbb Z$-valued random variables with law $\mu$. The reflected random walk $(X_n)$ is defined recursively by $X_0=x \in \mathbb N_0, X_{n+1}=\vert X_n+Y_{n+1}\vert$.  Under mild hypotheses on the law $\mu$, it is proved that,  for any $ y \in \mathbb N_0$, as $n \to +\infty$,  one gets $\mathbb P_x[X_n=y]\sim C_{x, y} R^{-n} n^{-3/2}$  when $\sum_{k\in \mathbb Z} k\mu(k) >0$  and 
$\mathbb P_x[X_n=y]\sim C_{ y} n^{-1/2}$  when $\sum_{k\in \mathbb Z} k\mu(k) =0$, for some constants $R, C_{x, y}$ and $C_y >0$.
 \end{abstract}


\section{Introduction}
 
We consider a sequence $(Y_i)_{i \geq 1}$ of $\mathbb Z$-valued independent and identically distributed random variables, with commun law $\mu$,   defined on a probability space $(\Omega, \mathcal F, \mathbb P)$. 

We denote by $(S_n)_{n \geq 0}$ the classical random walk  with law $\mu$ on $\mathbb Z$,  defined by $S_0=0$
and $S_n= Y_1+\cdots +Y_n$;  the canonical filtration associated with the sequence $(Y_i)_{i \geq 1}$ is denoted  $({\mathcal T}_n)_{n\geq 1}$. The {\bf reflected random walk} on $\mathbb N_0$ is defined by: for $X_0$ given and $\mathbb N_0$ valued, one sets
$$
\forall n \geq 0\quad X_{n+1}=\vert X_n+Y_{n+1}\vert .
$$
The process $(X_n)_{n \geq 0}$ is  a Markov chain on $\mathbb N_0$ with initial law $\mathcal L (X_0)$ and transition matrix $Q=(Q_{x, y})_{x, y \in  \mathbb N_0})$ given by 
$$
\forall x, y \geq 0 \quad q(x, y)=\left\{
\begin{array}{ll}
\mu(y-x)+\mu(y+x)& \mbox{\rm if}\quad y\neq 0\\
\mu(-x) &\mbox{\rm if}\quad y= 0
\end{array}
\right.
.
$$
When $X_0=x\ \mathbb P-$a.s., with $x \in \mathbb N_0$ fixed, the random walk $(X_n)_{n \geq 0}$ is denoted $(X_n^x)_{n \geq 0}$; the probability measure on $(\Omega, \mathcal T)$ conditioned to the event $[X_0=x]$ will be denoted $\mathbb P_x$ and the corresponding expectation $\mathbb E_x$.

We are interested with the behavior of the probabilities $\mathbb P_x[X_n=y], x, y \in \mathbb N_0$ as $n \to +\infty$; it is thus natural to consider the following generating function  $G$ associated with $(X_n)_{n \geq 0}$ and defined formaly as follows:
$$
\forall x, y \in \mathbb N_0, \forall s \in \mathbb C\quad {\mathfrak G}(s\vert x, y):= \sum_{n\geq 0}\mathbb P_x[X_n=y]s^n.
$$
The radius of convergence  $R$ of this series is of course $\geq 1$. 
  The reflected random walk is positive recurrent when 
$\mathbb E[\vert Y_i\vert ]<+\infty $ and $\displaystyle \mathbb E[ Y_i  ]<0$ (see \cite{PW} for instance and references therein) and consequently $R=1$;   it is also the case when 
 the $Y_i$ are centered, under the stronger assumption $\displaystyle \mathbb E[\vert Y_i\vert^{3/2} ]<+\infty $. A contrario,  when 
$\mathbb E[\vert Y_i\vert ]<+\infty $ et $\mathbb E[ Y_i  ] >0$,  as in the case of the classical random walk on $\mathbb Z$, it is natural to assume that $\mu$ has exponential moments$^($\footnote{namely we will assume that $\displaystyle \sum_{k \in \mathbb Z} r^k\mu(k) <+\infty$ for any $r>0$}$^)$ and, under this additional assumption, we will see that $R>1$.

The generating functions are of interest since we can often recover information  about the asymptotic behavior of probabilities, for instance by resorting a Tauberian theorem, e.g. that of Karamata; unfortunately, in this situation we have no way of obtaining the necessary information about these probabilities to apply such a Tauberian theorem (usually, we require that the sequence $p_n$ is monotone, which is far to be right in our situation) and we will employ the following theorem of Darboux: it requires more regularity of the generating function in a neighborhood of the singular point $z=R$ than does Karamata's theorem but no monotony type assumption: 

\begin{theorem}\label{darboux}
Let $\displaystyle \mathfrak G(s)= \sum_{n=0}^{+\infty} g_ns^n$ be a power series  with nonnegative coefficients $p_n$ and radius of convergence $R>0$. We assume that $\mathfrak G$ has no singularities in the closed disk $\Bigl\{s\in \mathbb C\slash \vert s\vert \leq R\Bigr\}$ except $s=R$ (in other words, $\mathfrak G$ has an analytic continuation to an open neighborhood of  the set $\Bigl\{s\in \mathbb C\slash \vert s\vert \leq R\Bigr\}\setminus \{R\}$) and that in a   neighborhood of $s=R$
\begin{equation}\label{singDARBOUX}
\mathfrak G(z) = \mathfrak{A}(s) (R-s)^{\alpha}+{\mathfrak B}(s)
\end{equation}
where $\mathfrak{A}$ and $\mathfrak{B}$ are analytic functions 
$^($\footnote{  in equation \ref{singDARBOUX}, this is the positive branch $s^\alpha$ which is meant, which implies that the branch cut is along the negative axis; so the branch cut for the function $\mathfrak G (s)$ is along the halfline $[R, +\infty]$}$^)$. Then 
\begin{equation}\label{asymDARBOUX}
g_n\sim { {\mathfrak A}(R) R^{1-n}\over \Gamma(-\alpha)n^{1+\alpha}} \quad {\it as} \quad n\to +\infty.
\end{equation}
\end{theorem}

This approach has been yet developed   by S. Lalley in the general context of {\it random walk with a finite reflecting zone}; the transitions $q(x, \cdot)$ of   Markov chains of this class are the ones of a classical random walk on $\mathbb N_0$ whenever $x \geq K$ for some $K\geq 0$. In our context of the reflected random walk on $\mathbb N_0$, it means that the support of $\mu$ is bounded from below (namely by $-K$); we will not assume this in the sequel and will thus not follow the same strategy than S. Lalley. The methods required for the analysis of random walks with non localized reflections are  more delicate, this is the aim of the present work for a particular such a  process.

The reflected random walk on $\mathbb N_0$ is characterized by the existence of reflection times. We have to consider the sequence $({\bf r}_k)_{k\geq 0}$  of successive reflection times; this is a sequence of waiting time with respect to the filtration $(\mathcal T_n)_{n \geq 0}$, defined by
$$
{\bf r}_0=0\quad {\rm and}\quad {\bf r}_{k+1}:= \inf\{n>{\bf r}_k : X_{{\bf r}_k}+Y_{{\bf r}_k+1}+\cdots+Y_n<0\} \quad{\rm for \ all} \quad k\geq 0.
$$
In the sequel we will often omit the index for ${\bf r}$ and denote the first reflection time ${\bf r}$.
If one assume 
$\mathbb E[\vert Y_i\vert ]<+\infty $ and $\mathbb E[ Y_i  ]\leq 0,$ one gets
 $\mathbb P_x[{\bf r}_k<+\infty]=1$ for all  $x \in \mathbb N_0$ et $k\geq 0$; on the contrary, when $\mathbb E[\vert Y_i\vert ]<+\infty $ and $\mathbb E[ Y_i  ]>0,$ one gets 
$\mathbb P_x[{\bf r}_k<+\infty]<1$ and in order to have $\mathbb P_x[{\bf r}_k<+\infty]>0$ it is necessary to assume that $\mu(\mathbb Z^{*-})>0$.

The following identity will be essential in this work, it can be stated in an elementary way : 
\begin{proposition} For all $x, y \in \mathbb N_0$ and  $s \in \mathbb C$, one gets 
\begin{equation}\label{factor1}
{\mathfrak G}(s\vert x, y)=\mathfrak E(s\vert x, y)+\sum_{w\in \mathbb N^*}\mathfrak R(s\vert x, w){\mathfrak G}(s\vert w, y),
\end{equation}
with

$\bullet$ for all $x, y \geq 0$   
$$\mathfrak E(s\vert x, y):=\sum_{n=0}^{+\infty} s^n\mathbb P_x[X_n=y, {\bf r}>n],
$$

$\bullet $ for all  $x  \geq 0$  and $w\geq 1$
\begin{eqnarray*}
\mathfrak R(s\vert x, w)&:=&\mathbb E_x[ 1_{[{\bf r}<+\infty, X_{{\bf r}}=w ]}s^{{\bf r}}]\\
&=& \sum_{n \geq 0} 
s^n\mathbb P[x+S_1\geq 0, \cdots, x+S_{n-1}\geq 0, x+S_n=-w].
\end{eqnarray*}
\end{proposition}
The generating function $\mathfrak E$ concerns the excursion of the Markov chain $(X_n)_{n \geq 0}$ before its first reflection and $\mathfrak R$  is related to the  process of reflection $(X_{{\bf r}_k})_{k \geq 0}$.

Proof. Let us decompose ${\mathfrak G}(s\vert x, y)$ into ${\mathfrak G}_1(s\vert x, y)+{\mathfrak G}_2(s\vert x, y)$ with
$$
{\mathfrak G}_1(s\vert x, y):=\mathbb E_x\left[ \sum_{n=0}^{{\bf r}-1}1_{\{y\}}(X_n) s^n\right] \quad \mbox{\rm and} \quad {\mathfrak G}_2(s\vert x, y):=\mathbb E_x\left[1_{[{\bf r}<+\infty]} \sum_{n={\bf r}}^{+\infty}1_{\{y\}}(X_n) s^n\right].
$$
One gets  ${\mathfrak G}_1(s\vert x, y)={\mathfrak E}(s\vert x, y)$ and, on the other side, by the strong Markov property, 
\begin{eqnarray*}
{\mathfrak G}_2(s\vert x, y) 
&=&  \sum_{k\geq 0}   \mathbb E_x\left[1_{[{\bf r}<+\infty]} 1_{[X_{{\bf r}+k}=y]} s^{{\bf r}+k}\right]\\
&=&  \sum_{k\geq 0} \sum_{w \in \mathbb N^*}  \mathbb E_x\left[1_{[{\bf r}<+\infty, X_{{\bf r}}=w]}s^{{\bf r}}
\mathbb P_w\left[X_k=y\right]s^k \right]\\
&=&  \sum_{w \in \mathbb N^*}  \mathbb E_x\left[1_{[{\bf r}<+\infty, X_{{\bf r}}=w]}s^{{\bf r}}\right]\times 
 \sum_{k\geq 0}\mathbb P_w\left[X_k=y\right]s^k \\
&=& \sum_{w\in \mathbb N^*}{\mathfrak R}(s\vert x, w){\mathfrak G}(s\vert w, y).
\end{eqnarray*}
\rightline{$\Box$}

By (\ref{factor1}), one easily sees that, to precise the asymptotic behavior of the $\mathbb P_x[X_n=y]$, it is necessary to control  the excursions of the walk between two successive reflection times. Note that this interrelationship among the  Green's functions $G, F$ and $H$ may be written as a single matrix equation involving matrix-valued generating functions.
For  $s\in \mathbb C$, let us denote ${\mathcal G}_s, {\mathcal E}_s$ and  ${\mathcal R}_s$  the following infinite matrices

$\bullet\quad 
{\mathcal G}_s=({\mathcal G}_s(x, y))_{x, y \in \mathbb N_0}$ with  ${\mathcal G}_s(x, y)={\mathfrak  G}(s\vert x, y)$ for all $x, y \in \mathbb N_0,$

$\bullet\quad 
{\mathcal E}_s=({\mathcal E}_s(x, y))_{x, y \in \mathbb N_0}$ with  ${\mathcal E}_s(x, y)={\mathfrak E}(s\vert x, y)$ for all $x, y \in \mathbb N_0,$

$\bullet\quad 
{\mathcal R}_s=({\mathcal R}_s(x, y))_{x \in \mathbb N_0, y\in \mathbb N^*}$ with  ${\mathcal R}_s(x, y)={\mathfrak R}(s\vert x, y).$

\noindent Thus for all $x, y \in \mathbb N_0$  and $s \in \mathcal C$, one gets 
\begin{equation} \label{factor2}
{\mathcal G}_s={\mathcal E}_s+{\mathcal R}_s{\mathcal G}_s.
\end{equation}
This shows that the Green functions ${\mathfrak G}(\cdot \vert x, y )$ may be computed when  $I-{\mathcal R}_s$ is invertible,  in which case one may write 
$$ {\mathcal G}_s=(I-{\mathcal R}_s)^{-1} {\mathcal E}_s.$$

Let us now introduce some general assumptions:

\noindent {\bf Hypotheses  H}: 

{\bf H1}: {\it  the measure $\mu$ is {\bf adapted} on $\mathbb Z$ (i-e the group generated by its support $S_\mu$ is equal to 
$\mathbb Z$)  and {\bf aperiodic} (i-e the group generated by $S_\mu-S_\mu$ is equal to 
$\mathbb Z$)

 {\bf H2}: the measure $\mu$  has exponential moment of any order (i.e. $\displaystyle \sum_{n\in \mathbb Z}r^n\mu(n) <+\infty$ for any $r \in ]0, +\infty[$)
  and $\displaystyle \sum_{n\in \mathbb Z} n \mu(n) \geq 0$.
  $^($\footnote{we can in fact consider weaker assumptions: there exist  $0<r_-<1<r^+$ such that $\displaystyle \hat{\mu}(r):= \sum_{n\in \mathbb Z}r^n\mu(n) <+\infty$ for any $r \in ]r_-, r_+[$ and $\mu_{r}$ reaches its minimum on this interval at a (unique) $r_0\in ]r_-, 1]$. We thus need  much more notations at the beginning, this complicates in fact  the understanding   of the proof and is not really of interest.}$^)$}

We now state the main result of this paper, which extends \cite{L} in our situation:

\begin{theorem}\label{MAIN}
Let $(Y_n)_{n \geq 1}$ be a sequence of $\mathbb Z$-valued independent and identically distributed random variables  with law $\mu$  defined on a probability space $(\Omega, \mathcal F, \mathbb P)$. Assume that $\mu$ satifies Hypotheses { \rm \bf H} and let $(X_n)_{n\geq 0}$ be the reflected random walk defined inductively by
$$
X_0=x \qquad \mbox{\rm and} \qquad  X_{n+1}=\vert X_n+Y_{n+1}\vert 
\quad \mbox{\rm for} \quad n \geq 0.
$$

$\bullet \quad$ If $\displaystyle \mathbb E[Y_1]=\sum_{k \in \mathbb Z} k\mu(k) =0$, then for any $ y \in \mathbb N_0$, there exists a constant $C_{  y}\in \mathbb R^{*+}$ such that, for any starting point $x \in \mathbb N_0$,  one gets 
$$\mathbb P_x[X_n=y]  \sim  {C_{y}\over \sqrt{n}}\quad {\rm as} \quad n\to +\infty,$$

$\bullet \quad$ If $\displaystyle \mathbb E[Y_1]=\sum_{k \in \mathbb Z} k\mu(k) >0$ then, for any $x, y \in \mathbb N_0$, there exists a constant $C_{x,  y}\in \mathbb R^{*+}$ such that 
$$
 \mathbb P_x[X_n=y]\sim  C_{x, y}{ \rho^n\over n^{3/2}}
$$ for some $\rho = \rho(\mu) \in ]0, 1[.$

  \end{theorem}
The constant $\rho(\mu)$ which appears in this statement  is the infimum over $\mathbb R$ of the generating function of $\mu$. We also know  the  exact value of the constants $C_y$ and $C_{x, y}, x, y \in \mathbb N_0$ which appear in the previous statement : see formulae (\ref{constantCy}) and (\ref{constantCx,y}).

\section{Decomposition of the trajectories and factorizations }\label{Decomposition of the trajectories}

In this section, we will consider the subprocess of reflections $(X_{{\bf r}_k})_{k \geq 0}$ in order to decompose the trajectories of the reflected random walk in several parts which can be analyzed.

We  first  introduce   some notations which appear classically in the fluctuation theory of $1$-dimensional random walks.  

\subsection{ On the fluctuations of a classical   random walk on $\mathbb Z$} 
Let $\tau^{*-}$  the first strict descending time of the random  walk $(S_n)_{n \geq 0}$:
$$
\tau^{*-}:= \inf\{n\geq 1/ S_n<0\} $$
(with\ the \ convention $  \inf \emptyset = +\infty$). The variable $\tau^{*-}$ is a stopping time with respect to the   filtration $(\mathcal T_n)_{n \geq 0}$.

We denote by $(T^{*-}_n)_{n \geq 0}$ the sequence of successive ladder descending epoch of the random walk $(S_n)_{n \geq 0}$ defined by
$T^{*-}_0=0$ and $T^{*-}_{n+1}=\inf \{k>T^{*-}_{n}/S_k< S_{T^{*-}_{n}}\}$ 
for $n\geq 0$. One gets in particular 
$T^{*-}_1=\tau^{*-}$; furthermore, setting $\tau^{*-}_n:= T^{*-}_{n}-T^{*-}_{n-1}$ for any $n \geq 1$, one may write $ T^{*-}_{n}= \tau^{*-}_1+\cdots+\tau^{*-}_n$ where $(\tau^{*-}_n)_{n \geq 1}$ is a sequence of  of independent and identically  random variables with law $\mu^{*-}:= \mathcal L(S_{\tau^{*-}})$.
 The potential associated to $\mu^{*-}$ is denoted by $U^{*-}$; one gets
$$U^{*-}(\cdot):= \sum_{n = 0}^{+\infty} \Bigl(\mu^{*-}\Bigr)^{*n}(\cdot)= \sum_{n = 0}^{+\infty}\mathbb E\Bigl[ \delta_{S_{T_n^{*-}}}(\cdot)\Bigr].$$

Similarly, we can  introduce   the first ascending time $\tau^{+}:= \inf\{n\geq 1/ S_n\geq 0\}$ of the random  walk $(S_n)_{n \geq 0}$ 
(with\ the \ convention $  \inf \emptyset = +\infty$) and  the 
 the sequence  $(T^{+}_n)_{n \geq 0}$ of successive ladder ascending epoch of  $(S_n)_{n \geq 0}$ defined by
$T^{+}_0=0$ and $T^{+}_{n+1}=\inf \{k>T^{+}_{n}/S_k\geq S_{T^{+}_{n}}\}$ 
for $n\geq 0$; as above,  one may write $ T^{+}_{n}= \tau^{+}_1+\cdots+\tau^{+}_n$ where $(\tau^{+}_n)_{n \geq 1}$ is a sequence of i.i.d.  random variables with law $\mu^+:=\mathcal L(S_{\tau^{+}})$. The potential associated with $\mu^+$ is denoted by $U^+$ ; one gets
$$U^{+}(\cdot):= \sum_{n = 0}^{+\infty} \Bigl(\mu^{+}\Bigr)^{*n}(\cdot)= \sum_{n = 0}^{+\infty}\mathbb E\Bigl[ \delta_{S_{T_n^{+}}}(\cdot)\Bigr].$$
.

We   need to control the  law of the couple $(\tau^{*-}, S_{\tau^{*-}})$ and thus introduce the ``characteristic''  function        $\varphi^{*-}$ defined formally by
$$\varphi^{*-} : (s, z)\mapsto \sum_{n \geq 1}s^n \mathbb E\left[1_{[\tau^{*-}=n]} z^{S_{n}}\right]\ $$  for $s, z \in \mathbb C$. In other words,  one gets
$$
\varphi^{*-}  (s, z)=\mathbb E[1_{[\tau^{*-}<+\infty]} s^{\tau^{*-}}z^{S_{\tau^{*-}}}];
$$
when the $Y_i$ are centered, we know that $\tau^{*-}$ is a.s. finite and the indicator function will be omitted in the sequel, otherwise we will modify suitably the choice of the law of the $Y_i$  and will pull back the study of $\varphi^{*-} $ in the centered case.

 We also introduce the characteristic function associated to the potential of $(\tau^{*-}, S_{\tau^{*-}})$, defined formally by
$$
\Phi^{*-}(s, z) =\sum_{k\geq 0}\mathbb E\left[s^{T_k^{*-}}z^{S_{T_k^{*-}}}\right]= \sum_{k\geq 0}
\varphi^{*-}(s, z)^k={1\over 1-\varphi^{*-}(s, z)}.
$$

There  is  be a natural duality between the open  half-line $\mathbb R ^{*-}$ and its complementary set $\mathbb R^+$; as above, we associate to the couple $(\tau^{+}, S_{\tau^{+}})$ the function $\varphi^+$ defined by
$$\varphi^{+}: (s, z)\mapsto \mathbb E[s^{\tau^{+}}z^{S_{\tau^{+}}}],$$   for $s, z \in \mathbb C$ with modulus $\leq 1$.  In fact, in a natural way will appear the ``potential''  associated with $(\tau^{+}, S_{\tau^{+}})$ and  whose 
``characteristic'' function $(s, z)\mapsto \Phi^+(s, z)$ is given by
$$
  \Phi^+(s, z):= \sum_{k\geq 0}\mathbb E\Bigl[s^{T_k^+} z^{S_{T_k^+}}\Bigr]
  =\sum_{k\geq 0} \varphi^+(s, z)^k={1\over 1-\varphi^+(s, z)}
$$
for complex numbers $s, z$ with modulus $<1$ (since in this case $\vert \varphi^+(s, z)\vert <1$). Notice that, by  a straightforward argument, called  {\it duality lemma} in the book by Feller \cite{F}, one also gets 
\begin{equation}
\Phi^{+}(s, z) =  \sum_{n \geq 0} s^n
\mathbb E\left[\tau^{*-}>n, z^{S_{n}}\right].
\end{equation}
 
We now introduce the corresponding generating functions $\mathfrak T^{*-}, \mathfrak U ^{*-}$ and   $\mathfrak U ^{+}$ defined by, for any $s \in \mathbb C, \vert s\vert \leq 1$ and $x \in \mathbb Z$
\begin{eqnarray*}
  \mathfrak T^{*-}( s \vert x) &=& \mathbb E\Bigl[ s^{\tau^{*-}}  1_{\{x\}}(S_{\tau^{*-} })\Bigr] 
=\sum_{n \geq 1}s^n{\mathbb P}\Bigl[  \tau^{*-}=n,  S_n=x\Bigr],
\\
 \mathfrak T^{+}( s \vert x) &=& \mathbb E\Bigl[ s^{\tau^{+}}  1_{\{x\}}(S_{\tau^{+} })\Bigr] 
=\sum_{n \geq 1}s^n{\mathbb P}\Bigl[  \tau^{+}=n,  S_n=x\Bigr],
\\
 \mathfrak U^{*-}( s \vert x) &=& \sum_{k\geq 0}\mathbb E\Bigl[s^{T_k^{*-}} 1_{\{x\}}(S_{T_k^{*-}})\Bigr]= \sum_{n \geq 0}s^n{\mathbb P}\Bigl[  \tau^{+}>n,  S_n=x\Bigr],
\\
 \mathfrak U^{+}( s\vert x)&=&\sum_{k\geq 0}\mathbb E\Bigl[s^{T_k^{+}} 1_{\{x\}}(S_{T_k^{+}})\Bigr]= \sum_{n \geq 0}s^n{\mathbb P}\Bigl[  \tau^{*-}>n,  S_n=x\Bigr].
\end{eqnarray*}
Note that $\mathfrak U^{*-}( s \vert x)  =0$ when $x \geq 0$ and 
$\mathfrak U^{+}( s\vert x)=0 $ when $x\leq -1$. 

We will first study the regularity of the Fourier transforms $\varphi^{*-}$ and $\varphi^+$ to describe the one of the functions $\mathfrak T^{*-}(\cdot \vert x)$ and $\mathfrak T^{+}(\cdot \vert x)$; to do this we will use the Wiener-Hopf factorization theory, in a quite strong version, in order to obtain some uniformity in the estimations we will need. We could adapt the same approach  for the functions $\mathfrak U^{*-}(\cdot \vert x)$ and $\mathfrak U^{+}(\cdot \vert x)$, but it is more difficult to control the behavior near $s=1$ of their respective Fourier transforms $\Phi^{*-}$ and $\Phi^+$. We will thus prefer to note that, for any $x \in \mathbb Z^{*-}$, the function $\mathfrak U^{*-}(\cdot \vert x)$ is equal to the finite sum 
$\displaystyle  \sum_{k=0}^{\vert x\vert}\mathbb E\Bigl[s^{T_k^{*-}} 1_{\{x\}}(S_{T_k^{*-}})\Bigr],$ since $T^{*-}_k\geq k$ a.s; the same remark does not hold for $\mathfrak U^{+}(\cdot \vert x)$  since $\mathbb P[S_{\tau^+}=0]>0$ but we will see that the  series  $\displaystyle  \sum_{k=0}^{+\infty}\mathbb E\Bigl[s^{T_k^{+}} 1_{\{x\}}(S_{T_k^{+}})\Bigr]$ converges exponentially fast and a similar approach will be developped.

It will be of interest to  consider the following square infinite matrices

 $\bullet \quad \mathcal  T^{*-}_s=\Bigl(\mathcal T^{*-}_s(x, y)\Bigr)_{x, y \in \mathbb Z^-}$ with $ \mathcal T^{*-}_s(x, y):= \mathfrak T^{*-}(s\vert y-x )$ for any $x, y \in \mathbb Z^-$,

 $\bullet \quad \mathcal U^{*-}_s=\Bigl(\mathcal U^{*-}_s(x, y)\Bigr)_{x, y \in \mathbb Z^-}$ with $\mathcal U^{*-}_s(x, y):= \mathfrak U^{*-}(s\vert y-x)$ for any $x, y \in \mathbb Z^-$.

The element of $\mathbb Z^-$ are labelled here in the decreasing order. Notice that the matrix $\mathcal  T^{*-}_s$ is strictly upper triangular; so  for any $x, y \in \mathbb Z^-$ one gets 
$\displaystyle \mathcal U^{*-}_s(x, y)=\sum_{k=0}^{\vert x-y\vert} (\mathcal  T^{*-}_s)^k(x, y)$.

 $\bullet \quad \mathcal  T^{+}_s=\Bigl(\mathcal T^{+}_s(x, y)\Bigr)_{x, y \in \mathbb N_0}$ with $ \mathcal T^{+}_s(x, y):= \mathfrak T^{+}(s\vert y-x )$ for any $x, y \in \mathbb N_0,$

 $\bullet \quad \mathcal U^{+}_s=\Bigl(\mathcal U^{+}_s(x, y)\Bigr)_{x, y \in \mathbb N_0}$ with $\mathcal U^{+}_s(x, y):= \mathfrak U^{+}(s\vert y-x)$ for any $x, y \in \mathbb N_0.$

 We will aso have $\displaystyle \mathcal U^{+}_s(x, y) = \sum_{k \geq 0}(\mathcal  T^{+}_s)^k(x, y)  $ for any $x, y \in \mathbb N_0$, the number of terms in the sum will not be finite in this case but it will not be difficult to derive the regularity of the function $s \mapsto \mathcal U^{+}_s(x, y) $ from the one of each term $s \mapsto   \mathcal  T^{+}_s(x, y)$.

In the sequel, we will consider the matrices $\mathcal T_s^{*-}$ and $\mathcal T_s^{+}$  as operators acting on $\left(\mathbb C^{\mathbb N_0}, \vert \cdot\vert_\infty\right)$; it  will  not be possible to  give sense to  the above  inversion formula on the Banach space of linear continuous operators acting on $\left(\mathbb C^{\mathbb N_0}, \vert \cdot\vert_\infty\right)$ and we will have to consider the action of these matrix and a larger space of $\mathbb C$-valued sequences.

In the following subsections,  we decompose both   the excursion of $(X_n)_{n \geq 0}$ before the first reflection and   the process of reflections $(X_{\bf r_k})_{k \geq 0}$ in terms of  quantities introduced here.

\subsection{The approach process and the matrices $\mathcal T_s$}

The trajectories of the reflected random walk are governed by the strict descending ladder epoch of the corresponding  classical random walk on $\mathbb Z$, and the generating function
$\mathfrak T^{*-}$ introduced in the previous section will be essential in the sequel. Since the staring point may be any $x \in \mathbb N_0$, we have to consider the first time at which the random walk $(X_n)_{n\geq 0}$ goes on the ''left'' on the initial point (with eventually a reflexion at this time, in which case the arrival point may be $>x$), that is the strict descending ladder epoch $\tau^{*-}$ of the random walk $(S_n)_{n \geq 0}$. We thus introduce the   matrices $\mathcal T_s$  which contains a lot of information for the reflected random walk, defined by $\mathcal T_s=\Bigl(\mathcal T_s(x, y)\Bigr)_{x, y \in \mathbb N_0}$ with 
\begin{equation}\label{approach}
\forall x, y \in \mathbb N_0\qquad \mathcal T_s(x, y):= \mathfrak T^{*-}(s\vert y-x).
\end{equation}
Notice that the matrices $\mathcal T_s$ are  strictly lower triangular.

\subsection{The excursion before the first reflection}
Recall that  the function $\mathfrak E$ is defined by
$$ \forall x, y \in \mathbb N_0, \forall s \in \mathbb C \qquad 
\mathfrak E(s\vert x, y):= \sum_{n\geq 0} s^n\mathbb P_x\left[ {\bf r}>n, X_n=y\right].
$$
We have the following identity: for all $s \in \mathbb C$ and $x, y \in \mathbb N_0$
\begin{equation*} 
\mathfrak E (s\vert x, y)=\mathfrak U^{+}( s\vert y-x )+\sum_{w=0}^{x-1}
\mathfrak T^{*-}(s\vert w-x)\mathfrak E(s\vert w, y ).
\end{equation*}
As above, we introduce the square infinite matrices    $\mathcal E_s=\Bigl(\mathcal E_s(x, y)\Bigr)_{x, y \in \mathbb N_0}$, 
with  $\mathcal E_s(x, y):= \mathfrak E(s\vert x, y)$ for any $x, y \in \mathbb N_0$,   and   rewrite this identity as follows
$$
\mathcal E_s=\mathcal U^+_s + \mathcal  T_s \mathcal E_s.
$$
Since $\mathcal  T_s$ is strictly lower triangular, the matrix  $I-\mathcal  T_s$ will be  invertible (in a suitable space to be precised) and one will get
\begin{equation}\label{factor1bis}
\mathcal E_s= \Bigl(I-\mathcal  T_s\Bigr)^{-1}\mathcal U^+_s.
\end{equation}
In the follwing sections, we will give sense to this inversion formula and describe the regularity in $s$ of the matrix-valued function $s \mapsto \mathcal E_s$.

 \subsection{The process of reflections}
 
 Under the hypothesis $\mathbb P[\tau^{*-} <+\infty]=1\ ^(\footnote{this condition is satisfied for instance when $\mathbb E[\vert Y_i\vert ]<+\infty $ and  $\mathbb E[ Y_i  ]\leq 0$.}^)$, the distribution law of the variable $S_{\tau^{*-}}$ is denoted $\mu^{*-}$ and its potential    $\displaystyle U^{*-}:= \sum_{n \geq 0} (\mu^{*-})^{*n}$; all the waiting times $T_n^{*-}$ are thus a.s. finite and one gets  $(\mu^{*-})^{*n}={\mathcal L}(S_{T_n^{*-}})$, furthermore, for any $x \in \mathbb N_0$ the successive reflection times  $ {\bf r}_k, k \geq 0, $ are also a.s. finite.  The process   $(X_{{\bf r}_k})_{k \geq 0}$ appears in a crucial way in \cite{PW}  to study the recurrence/transience properties of the reflected walk; indeed, we have the
 \begin{fact} \label{reflectedPW}\cite{PW}
Under the hypothesis  $\mathbb P[\tau^{*-} <+\infty]=1$, the process of reflections $(X_{{\bf r}_k})_{k \geq 0}$ is a Markov chain on  $\mathbb N_0$ with transition probability   ${\mathcal R}$ given by
\begin{equation}\label{transitionreflection}
\forall x   \in \mathbb N_0,  \forall y \in \mathbb N_0  \qquad 
{\mathcal R}(x, y)= \left\{
\begin{array}{lll }
0 & {\rm if} & y=0\\
\displaystyle \sum_{0}^{x} U^{*-}(-w)\mu^{*-}( w-x-y) & {\rm if} & y\geq 1
\end{array}
\right.
\end{equation}
Furthermore, the measure  $\nu_{\bf r}$  on $\mathbb N^*$ defined by
\begin{equation}\label{stationary}
\forall x \in \mathbb N^*\quad 
\nu_{\bf r}(x):= \sum_{y=1}^{+\infty}\left(
{\mu^{*-}(-x)\over 2}+\mu^{*-}\Bigl(]-x-y, -x[\Bigr)+{\mu^{*-}(-x-y)\over 2}\right)\mu^{*-}(-y)
\end{equation}
is stationary for $(X_{{\bf r}_k})_{k \geq 0}$ and is unique up to a multiplicative constant; it is finite  as soon as  $\displaystyle \mathbb E[\vert S_{\tau^{*-}}\vert ] = \sum_{k \geq 1} k\mu^{*-}(-k)<+\infty$.
\end{fact}
This statement is a bit different from the one in \cite{PW} since we assume here that at the reflection time the process $(X_n)_{n \geq 0}$ belongs to $\mathbb N^*$; nevertheless, the proof goes exactly along the same lines.
This result is crucial in the sequel in order to control the spectrum of the stochastic infinite matrix $ \mathcal R=\Bigl(\mathcal R(x, y)\mathcal )_{x, y \in \mathbb N_0}$; namely, we have the

\begin{property}
\label{quasicompact} There exists a constant $\kappa\in ]0, 1[$ such that, for any $x\in \mathbb N_0$ and $y \in \mathbb N^*$ one gets
$${\mathcal R} (x, y) \geq \kappa \mu^{*-}(-y).$$
In particular, the operator  $\mathcal R$  acting on $(\mathbb C^{\mathbb N_0}, \vert \cdot\vert_\infty)$ is quasi-compact : more precisely,  the eigenvalue $1$ is   simple,  with associated eigenvector $h=(1)_{n \in \mathbb N_0}$ and the rest of the spectrum is included in a disk of radius $\leq 1-\kappa$.

Furthermore, for any $K>1$, the operator $\mathcal R$ acts also on the Banach space $(\mathbb C^{\mathbb N_0}, \vert \cdot\vert_K)$, where $\vert \cdot \vert _K$ denotes the  norm defined by
\begin{equation}\label{normK}
\forall {\bf a}=(a_x)_{x \in \mathbb N_0} \in \mathbb C^{\mathbb N_0}\qquad 
\vert {\bf a}\vert_K:= \sup_{x\in \mathbb N_0} {\vert a_x\vert \over K^x}, 
\end{equation}
the eigenvalue $1$ is simple with associated eigenvector $h$ and the rest of the spectrum  of $\mathcal R$ acting on $(\mathbb C^{\mathbb N_0}, \vert \cdot\vert_K)$ is included in a disk of radius $\leq 1-\kappa$.

\end{property}
Proof. Let $N_\mu:= \inf\{k\leq -1\slash \mu\{k\}>0\}$ (with $N=-\infty$ is the support of $\mu$ is not bounded from below).
Since $\mu$ is adapted, one gets $\mu^{*-}(k)>0$ for any $k\in \{-N_\mu, \cdots, -1\}$ (and any $k \in \mathbb Z^{*-}$ when $N_\mu=-\infty$); as a direct consequence, one gets $U^{*-}(k) >0 $ for any $k \in \mathbb Z^{*-}$.
In fact, by the 1-dimensional renewal theorem, one knows that
$\displaystyle \lim_{k\to -\infty}U^{*-}(k)={1\over - \mathbb E [S_{\tau^{*-}}]}>0$ since $\mathbb E [S_{\tau^{*-}}]>-\infty$ when $\mu$ has exponential moments; it readily follows that
$\displaystyle 
\kappa:= \inf_{k \in \mathbb Z^{-}} U^{*-}(k) >0.
$
Using (\ref{transitionreflection}), one may thus write, for any $x \in \mathbb N_0$ and $y \in \mathbb N^*$
$$
\mathcal R(x, y)\geq U^{*-}(x)\mu^{*-}(-y)\geq \kappa \mu^{*-}(-y).
$$
The matrix $(\mathcal R(x, y)_{x, y \in \mathbb N_0}$ thus satisfies the so-called ``Doeblin condition'' and it is quasi-compact on $(\mathbb C^\mathbb N_0, \vert\cdot\vert _\infty)$ (see for instance \cite{B} for a precise statement).

The same spectral property holds on $(\mathbb C^{\mathbb N_0}, \vert \cdot\vert_K)$ since $\mu^{*-}$ has exponential moment of any order, which allows to check  
 by a straightforward computation  that 
 $$
\sup_{x \in \mathbb N_0}\sum_{y\in \mathbb N_0} {\mathcal R}(x, y) K^y<+\infty. 
$$
\rightline{ $\Box$} 

{\bf For technical reasons which will appear in Section 4}, we will replace the function $x\mapsto  K^x$ by a function denoted also $K$ which satisfies the  following conditions
\begin{equation}\label{condition K}
\forall x \in \mathbb N_0  \quad K(x) \geq 1, \quad\mathcal R K(x)\leq 1 \quad {\rm and}  \quad K(x) \sim K^x.
\end{equation}
 It suffices to consider the  function 
$x \mapsto \Bigl(1\vee {K(x)\over M}\Bigr)$ with $\displaystyle M:= \sup_{x\in \mathbb N_0} \sum_{y \in \mathbb N^*}\mathcal R(x, y) K^y$ (we now that $M<+\infty$ by   proof of Property \ref{quasicompact}). {\bf The set of fonctions which satisfy the conditions (\ref{condition K}) will be denoted $\mathcal K(K)$.}

We now explicit the connection between  $\mathcal R_s$ and the matrix $\mathcal  T_s$ introduced above; namely, there exists a  similar factorization identity than (\ref{factor1}) for the process of reflection. Using the fact that the first reflection time may appear or not at  time $\tau^{*-}$, one may write: for all $s \in \mathbb C$ and $x \in \mathbb N_0$ and $y \in \mathbb N^*$
\begin{equation}\label{factor2}
\mathfrak R (s\vert x, y)=\mathfrak T( s\vert -x-y )+\sum_{w=0}^{x-1}
\mathfrak T(s\vert w-x)\mathfrak R(s\vert w, y ),
\end{equation}
which leads to the following equality:
\begin{equation}\label{factor2bis}
\mathcal R_s= \Bigl(I-\mathcal  T_s\Bigr)^{-1}\mathcal V_s
\end{equation}
where we have set $\mathcal V_s=\Bigl(\mathcal V_s(x, y)\Bigr)_{x, y\in \mathbb N_0}$
with 
\begin{equation}\label{mathcalVs}
\mathcal V_s(x, y):= \left\{
\begin{array}{ll}
0 &{\rm if} \ y=0\\
\mathfrak T^{*-}( s\vert -x-y ) &{\rm if} \ y\in \mathbb N^*.
\end{array}
\right.
\end{equation}

The crucial point in the sequel will be thus to describe the regularity of the maps $s \mapsto \mathcal T_s, s \mapsto \mathcal  V_s$ and $ s \mapsto \mathcal U^+_s$ near the point $s=1$. 
We will first detail the centered case; the main ingredient is the classical Wiener-Hopf factorization which permits to control both functions 
$\varphi^{*-} $ and $\varphi^{+}$.

Another essential point will be to describe the one of the maps  $ \left(I-\mathcal  T_s\right)^{-1}$ and $ \left(I-\mathcal R_s\right)^{-1}$ and this question is related to the description of the spectrum of the operators $\mathcal  T_s$ and $\mathcal R_s$ when $s$ is closed to $1$: this is not difficult for $\mathcal  T_s$ since it is a strictly lower triangular matrix but more delicate for $\mathcal R_s$ in the centered case where $\mathcal R=\mathcal R_1$ is a Markov operator.

\section{A strong version of the Wiener-Hopf factorization and its applications to classical random walks}

 \subsection{Introduction and notations}The Wiener-Hopf factorization proposes a decomposition of the space-time characteristic function $(s, z) \mapsto 1-s\mathbb E[z^{Y_1}]= 1-s\hat{\mu}(z)$ in terms of $\varphi^{*-}$ and $\varphi^{+}$; namely,  one gets; for all 
  $s, z \in \mathbb C$ with modulus $<1$
  \begin{equation}\label{WH1}
  1-s \hat{\mu}(z)= \Bigl(1-\varphi^{*-} (s, z)\Bigr) \Bigl(1-\varphi^{+} (s, z)\Bigr).
\end{equation}
In \cite{EP}, we already use this factorization in order to state local limit theorems for fluctuations of the random walk $(S_n)_{n \geq 0}$; we first propose another such a decomposition, and, by identification of the corresponding factors, we obtain another expression  for each of the functions $\varphi^{*-}$ and $\varphi^{+}$.This new expression allows us to use  elementary arguments coming from entire functions theory in order to describe for instance the  asymptotic behavior of the sequences 
$\Bigl(\mathbb P[S_n=x, \tau^{*-}=n]\Bigr)_{n \geq 1}$ and 
$\Bigl(\mathbb P[S_n=y, \tau^{*-}>n]\Bigr)_{n \geq 1}$ for any $x \in \mathbb Z^{*-}$ and $y \in \mathbb Z^+$. 

In the present situation, we need first to obtain similar results than in \cite{EP}  but in terms of regularity  with respect to the variable $s$ of the functions   $\varphi^{*-}$ and $\varphi^{+}$ around the unit circle, with a precise description of their singularity near the point $s=1$; by the identity (\ref{factor1})  we will show that these properties spread to the function $G(s\vert x, y)$, which allows us to conclude, using the classical Darboux's method for entire functions.

 We will assume that  the law $\mu$ as exponential moment of any order, i.e. $\displaystyle \sum_{n \in \mathbb Z}r^n \mu(n) <+\infty$ for any $r \in \mathbb R^{*+}$; it readily implies that its generating function $\hat{\mu}: z\mapsto \displaystyle \sum_{n \in \mathbb Z}z^n \mu(n) $ is analytic on $\mathbb C^*$; furthermore, its restriction to $]0, +\infty[$ is  strictly convex and one gets $\displaystyle \lim_{r\to +\infty} \hat{\mu}(r)=\lim_{r\to <0} \hat{\mu}(r)=+\infty$ as soon as  $\mu $ charges $\mathbb Z^{*+}$ and $\mathbb Z^{*-}$. In particular, under these conditions, there exists a unique $r_0 >0$ such that $\displaystyle \hat{\mu}(r_0)=\inf_{r>0} \hat{\mu}(r)$; one gets $\hat{\mu}'(r_0)=0, \hat{\mu}''(r_0)>0$ and sets $\rho_0:= \hat{\mu}(r_0)$. Note that $\rho_0=1$ when $\mu$ is centered and $\rho_0\in ]0, 1[$ otherwise; we will set $R_\circ:= {1\over \rho_0}$.

We now fix $0<r_-<r_0<r_+<+\infty$ and will denote by  ${\bf L}= {\bf L}[r_-, r_+]$ the space of functions $F: \mathbb C^*\to \mathbb C$ of the form 
$F(z):= \sum_{n \in \mathbb Z} a_nz^n$ for some (bilateral)-sequence $(a_n)_{n \in \mathbb Z}$ such that  $\displaystyle   \sum_{n \leq 0} \vert a_n\vert  r_-^n+\sum_{n \geq 0} \vert a_n\vert r_+^n<+\infty$; the elements of  $\bf L$ are   called {\bf Laurent functions} on the annulus $[r_-, r_+]:= \{r_-\leq \vert z\vert\leq  r_+\}$ and  the Banach space  $({\bf L}, \vert\cdot\vert_\infty)$ $^($\footnote{where $\vert\cdot\vert_\infty$ denotes the  norm  of uniform convergence on the annulus $\{r_-\leq \vert z\vert\leq  r_+\}$}$^)$   contains the function $\hat{\mu}$ defined above. 

\subsection{The centered case}
Lets us first consider the centered case: $\mathbb E[Y_i]=\hat{\mu}'(1)=0$; we thus have $r_0=1$ and $\rho_0=R_\circ=1$. 
Under the aperiodicity condition on $\mu$, one gets $\vert 1-s \hat{\mu}(z)\vert >0$ for any $z \in \mathbb C^*, \vert z\vert =1$, and $s$ such that $\vert s\vert \leq 1$,  excepted $s=1$; it follows that for any $z \in \mathbb C^*, \vert z\vert =1$, the function $\displaystyle s\mapsto {1\over 1-s \hat{\mu}(z)}$ may be analytically  extended on the  set $\{s\in \mathbb C/ \vert s\vert \leq 1+\delta\}\setminus[1, 1+\delta[$ for some $\delta>0$. 
On the other hand, setting  $\sigma^2:= \mathbb E[Y_i^2]$, one gets $\hat{\mu}''(1)=\sigma^2>0$.  One thus gets, setting $\Psi(s,z):=1-s\hat{\mu}(z),$
$$
{\partial\over \partial z}\Psi(1,1)=0\quad{\rm and} \quad {\partial^2\over \partial z^2}\Psi(1,1)=\sigma^2>0.
$$
The Weiertrass preparation lemma thus implies that, on a neighborhood of $(1, 1)$ one may write 
$$
1-s\hat{\mu}(z)=  \left( z^2+b(s)  z+c(s)\right)\Psi(s,z)
$$
with $\Psi$ analytic on $\mathbb C\times \mathbb C^*$ and $\Psi\neq 0$ on a neighborhood of $(1, 1)$. One gets   
$$z^2+b(s)z+c(s)= \left(z-z_-(s))(z-z_+(s)\right),$$ with $ z_-(s)<1<z_+(s)$ when $s\in [0, 1[$ and $z_-(1)=z_+(1)=1.$ 

In order to solve this last equation, we fix  the  principal determination   of the function $Z \mapsto \sqrt{Z}$ $^($\footnote{for $Z$ in $\mathbb C\setminus \mathbb R^{*-}$, writing $Z= \vert Z\vert e^{i\theta}$ for some $-\pi<\theta< \pi$, one set $\sqrt{Z}= \sqrt{\vert Z\vert}
e^{i\theta/2}$}$^)$ in such a way $s \mapsto \sqrt{1-s}$ is well defined on the set ${\mathcal O}_\delta(1):= B(1, \delta)\setminus[1, 1+\delta[$. It follows that the functions $z_{\pm} $ admit the  analytic expansion  $\displaystyle z_{\pm}(s)=1+\sum_{n \geq 1}(\pm1)^n \alpha_n(1-s)^{n/2}$  on ${\mathcal O}_\delta(1)$  and the equality $\displaystyle \hat{\mu}(z_{\pm}(s))={1\over s}=\sum_{n \geq 0}(1-s)^n$ valid for $0\leq s<1$
leads to  $\displaystyle \alpha_1= { \sqrt{2}\over \sigma}$.

This type of singularity of the functions $z_\pm$  near $s=1$ is essential in the sequel because it contains the one of the functions $\varphi^{*-} (s, z)$ and $\varphi^{+} (s, z)$  near $(1, 1)$. The Wiener-Hopf factorization   has several versions in the literature; we emphasize here that we need some kind of uniformity with respect to the parameter $z$ in the local expansion of the function $\varphi^{*-} $ near $s = 1$, this is why we consider the map $s\mapsto  \varphi^{*-} (s, \cdot)$ with values in  ${\bf L}[r_-, r_+]$.  It is proved in particular in \cite{B} (see also \cite{Pr} for a more precise statement, in the context of Markov walks) that there exists $\delta>0$ such  that the function 
$\displaystyle 
s \mapsto \Bigl( z\mapsto \phi^{*-}(s,z):=  {1-\varphi^{*-}(s, z)\over z-z_-(s)}\Bigr)
$
is analytic on the open ball $B(1, \delta)\subset \mathbb C$, with values in $L[r_-, r_+]$.
Setting $\displaystyle \phi^{*-}(s, \cdot)= \sum_{k\geq 0}\phi^{*-}_{(k)}(1-s)^k$ for $\vert 1-s\vert <\delta$ and $\phi^{*-}_{(k)}\in {\bf L}[r_-, r_+]$ and using the local expansion $z_-(s)=1-{{\sqrt 2}\over \sigma}\sqrt{1-s}+\cdots $, one thus gets  for $\delta$ small enough  and $s \in {\mathcal O}_\delta(1)$
$$
\varphi^{*-}(s, \cdot ) =   \varphi^{*-}(1, \cdot)+ \sum_{k\geq 1}\varphi^{*-}_{(k)}(1-s)^{k/2}$$
with $\displaystyle \sum_{k\geq 0}\vert \varphi^{*-}_{(k)}\vert_\infty \delta^k<+\infty$ and 
$\displaystyle \varphi^{*-}_{(1)} : z \mapsto {\sqrt{2}\over \sigma}\times {1-\mathbb E[z^{S_{\tau^{*-}}}]\over 1-z} $.

We summarize the informations we will need  in the following

\begin{proposition} \label{varphicentered} For any $r_-<1<r_+$, the function $s\mapsto \varphi^{*-}(s, \cdot)$  has an  analytic continuation to an open neighborhood of $\overline{B(0, 1)}\setminus\{1\}$   with values in ${\bf L}[r_-, r_+]$; furthermore, for  $\delta>0$      , this function is analytic in the variable $\sqrt{1-s}$ 
 on   the set  $ {\mathcal O}_\delta(1)=B(1, \delta)\setminus[1, 1+\delta[$ and  its local expansion of order 1 in ${\bf L}[r_-, r_+]$  is given by 
\begin{equation}\label{local pour varphi}
\varphi^{*-}(s, \cdot ) =  \varphi^{*-}(1, \cdot) +  \sqrt{1-s}\  \varphi^{*-}_{(1)}(\cdot)+{\bf O}(s, \cdot ) \end{equation}
with 
$\quad \displaystyle \varphi^{*-}_{(1)} : z \mapsto {\sqrt{2}\over \sigma}\times {1-\mathbb E[z^{S_{\tau^{*-}}}]\over 1-z} $
and ${\bf O}(s, \cdot )$ uniformly bounded in ${\bf L}[r_-, r_+]$.
\end{proposition}

A similar statement  holds for the function $\varphi^{+}$; in particular, the local expansion  near $s=1$ follows from the one of the root $z_+(s)$, namely $z_+(s)=1+{\sqrt{2} \over \sigma}\sqrt{1-s}+\cdots$. We may thus state the

\begin{proposition} \label{varphicentered} For any $r_-<1<r_+$, the function $s\mapsto \varphi^{+}(s, \cdot)$  has an  analytic continuation to an open neighborhood of $\overline{B(0, 1)}\setminus\{1\}$   with values in ${\bf L}[r_-, r_+]$; furthermore, for  $\delta>0$ small enough, this function is analytic in the variable $\sqrt{1-s}$ 
 on   the set  $ {\mathcal O}_\delta(1)=B(1, \delta)\setminus[1, 1+\delta[$ and  its local expansion of order 1 in ${\bf L}[r_-, r_+]$  is given by 
\begin{equation}\label{local pour varphi}
\varphi^{+}(s, \cdot ) =  \varphi^{+}(1, \cdot) +  \sqrt{1-s}\  \varphi^{+}_{(1)}(\cdot)+{\bf O}(s, \cdot ) \end{equation}
with 
$\quad \displaystyle \varphi^{+}_{(1)} : z \mapsto -  {\sqrt{2}\over \sigma}\times {1-\mathbb E[z^{S_{\tau^{+}}}]\over 1-z} $
and ${\bf O}(s, \cdot )$ uniformly bounded in ${\bf L}[r_-, r_+]$.
\end{proposition}

\subsection{The maps $s\mapsto \mathfrak T^{*-}( s \vert x)$ and $ s \mapsto \mathfrak T^+( s\vert x)$ for $x \in \mathbb Z$}
We use here the   inverse Fourier's formula: for any $x\in \mathbb Z^{*-}$ and $s\in \mathbb C, \vert s\vert <1$, one gets, by a Fubini type argument,
\begin{eqnarray*}
\mathfrak T^{*-}( s \vert x)&=& \mathbb E\left[s^{\tau^{*-}}1_{\{x\}}(S_{\tau^{*-}})\right]
\\
&=&  \mathbb E\left[s^{\tau^{*-}}{1\over 2i\pi}\int_{\mathbb T}z^{S_{\tau^{*-}}-x-1}dz\right]
\\
&=& {1\over 2i\pi}  \int_{\mathbb T}z^{-x-1}\varphi^{*-}(s, z) dz.
\end{eqnarray*}
Similarly   
$\mathfrak T^{+}( s \vert x)= {1\over 2i\pi} \int_{\mathbb T}z^{-x-1}\varphi^{+}(s, z) dz $ for any $x \in \mathbb N_0$. We will    apply Propositions  \ref{varphicentered} and \ref{varphicentered} and first identify the coefficients which appears in the local expansion as Fourier transforms of some known measures;  let us denote

$\bullet\ \delta_x$   the Dirac mass at $x \in \mathbb Z$, 

 $\bullet\ \lambda^{*-}=\displaystyle \sum_{x\leq -1} \delta_x$  the counting measures  on $\mathbb Z^{*-}$ 
 
 $\bullet\ \lambda^+=\displaystyle \sum_{n\geq 0} \delta_x$   the counting measures  on $\mathbb N_0$.
 
 One easily checks that $\displaystyle z \mapsto  {1-\mathbb E[z^{S_{\tau^{*-}}}]\over z-1}$ and  $   z \mapsto  {1-\mathbb E[z^{S_{\tau^{+}}}]\over  1-z}$    are  the generating functions  associated respectively  with the measures $
\displaystyle (\delta_0- \mu^{*-})* \lambda^{*-}$ and $
\displaystyle (\delta_0- \mu^+)* \lambda^{+};$ we may thus  state the following
\begin{proposition}\label{localexpansion}
There exists an open neighborhood $\Omega$ of $\overline{B(0, 1)}\setminus\{1\}$ such that, for any $x \in \mathbb Z $, the functions  $s\mapsto \mathfrak T^{*-}( s \vert x):=\mathbb E[s^{\tau^{*-}}1_{\{x\}}(S_{\tau^{*-}})]$  and 
$s\mapsto \mathfrak T^{+}( s \vert x):=\mathbb E[s^{\tau^{+}}1_{\{x\}}(S_{\tau^{+}})]$ 
have an  analytic continuation to   $\Omega$; furthermore, for  $\delta>0$ small enough, these functions are analytic in the variable $\sqrt{1-s}$ 
 on   the set  $ {\mathcal O}_\delta(1)$    and   their  local expansions of order 1  are given by 
\begin{equation}\label{mathfrakTx}
\mathfrak T^{*-}( s \vert x) =  \mu^{*-}(x) - \sqrt{1-s}{\sqrt{2}\over \sigma} \mu^{*-}\Bigl(]-\infty, x]\Bigr) +(1-s) \ {\bf O}(s\vert x)
\end{equation}
and 
\begin{equation}\label{local pour mathfrakT(x)}
\mathfrak T^{+}( s \vert x) =  \mu^{+}(x) - \sqrt{1-s}{\sqrt{2}\over \sigma} \mu^{+}\Bigl(]x, +\infty[\Bigr) +(1-s) \ {\bf O}(s\vert x)
\end{equation}
with 
 ${\bf O}(s\vert x)$ analytic in the variable $\sqrt{1-s}$ and uniformly bounded in $s \in {\mathcal O}_\delta(1)$ and $x \in \mathbb Z$.

Furthermore, for any $K>1$, there exists a constant ${\bf O}>0$ such  that 
\begin{equation}\label{uniformO}
K^{\vert x\vert} \Big\vert \mathfrak T^{*-}( s \vert x)\Big\vert \leq {\bf O}, \quad 
K^{\vert x\vert} \Big\vert   \mathfrak T^{+}( s \vert x)\Big\vert \leq {\bf O} \quad {\rm and}\quad K^{\vert x\vert} \Big\vert { \bf O}(s\vert x)\Big\vert  \leq {\bf O}.
\end{equation}
for any $s \in \Omega \cup {\mathcal O}_\delta(1)$ and $x \in \mathbb Z$.

 \end{proposition} 
Proof. The analyticity property and the local expansions (\ref{mathfrakTx}) and (\ref{local pour mathfrakT(x)}) are direct consequences of Propositions  \ref{varphicentered}  and \ref{varphicentered}. To establish  for instance the first inequality in 
(\ref{uniformO}), we use  the fact that  for $s \in \Omega \cup {\mathcal O}_\delta(1)$, the function $z \mapsto \varphi^{*-}(s, z)$ is analytic  on any annulus $\{z\in \mathbb C/ r_-< \vert z\vert < r_+\}$ with $0<r_- <1<r^+$ and so, for any $K>1$ and $x \in \mathbb Z^{*-}$, one gets 
$$
\mathfrak T^{*-}( s \vert x)={1\over 2i\pi}\int_{\mathbb T}z^{-x-1}\varphi^{*-}(s, z) dz={1\over 2i\pi}\int_{\{z/\vert z\vert =1/K\}}z^{-x-1}\varphi^{*-}(s, z) dz.
$$
So $\displaystyle \Big\vert \mathfrak T^{*-}( s \vert x)\Big\vert\leq  {K^{-\vert x\vert-1}\over 2\pi}\times  \sup_{\stackrel{s \in \Omega \cup {\mathcal O}_\delta(1)}{\vert z\vert =1/K}} \vert \varphi{*-}(s, z)\vert.
$ The same argument holds for the quantities $ \mathfrak T^{+}( s \vert x)$ and ${ \bf O}(s\vert x)$.

\rightline{$\Box$}

\subsection{ The  coefficient maps $s \mapsto \mathcal T_s^{*-}(x, y)$ and $s \mapsto \mathcal T_s^+(x, y)$ for $x, y \in \mathbb Z$}

We  first analyze here  the consequences of the previous statement for   the matrices coefficients $\mathcal T_s^{*-}(x, y)$ and $\mathcal T_s^+(x, y)$. 
We have the 
\begin{proposition} \label{Ts(x,y)}
There exists an open neighborhood $\Omega$ of $\overline{B(0, 1)}\setminus\{1\}$ such that  for any $x, y \in \mathbb Z$, the functions  $s\mapsto \mathcal T^{*-}_s (x, y)$  and 
$s\mapsto \mathcal T_s^{+}(x, y)$ 
have an  analytic continuation to   $\Omega$; furthermore, for  $\delta>0$ small enough, these functions are analytic in the variable $\sqrt{1-s}$ 
 on   the set  $ {\mathcal O}_\delta(1)$    and   their  local expansions of order 1  are given by 
\begin{equation}\label{local pour T(x,y)}
\mathcal T^{*-}_s(x,y) =  {\mathcal T}^{*-}(x,y) + \sqrt{1-s}\ \widetilde{\mathcal T}^{*-}(x,y) +(1-s) \ {\bf O}_s(x,y)
\end{equation}
and 
\begin{equation}\label{local pour Tplus}
\mathcal T^{+}_s(x,y) =  {\mathcal T}^{+}(x,y) + \sqrt{1-s}\ \widetilde{\mathcal T}^{+}(x,y) +(1-s) \ {\bf O}_s(x,y)
\end{equation}
where

$\bullet\quad  \displaystyle \mathcal T^{*-}(x, y)= \mu^{*-}(y-x),$

\vspace{2mm}

$\bullet\quad    \displaystyle 
  \widetilde{\mathcal T}^{*-}(x, y)= -{\sqrt{2}\over \sigma}\mu^{*-}\Bigl(]-\infty, y-x]\Bigr), $

\vspace{2mm}

$\bullet\quad \mathcal T^{+}(x, y)= \mu^{+}(y-x),$

\vspace{2mm}

$\bullet\quad  \displaystyle \widetilde{\mathcal T}^{+}(x, y)= -{\sqrt{2}\over \sigma}\mu^{+}\Bigl(]y-x, +\infty[\Bigr), $

\vspace{2mm}

$\bullet\quad $   ${\bf O}_s(x, y)$ is analytic in the variable $\sqrt{1-s}$  for $s \in {\mathcal O}_\delta(1)$.

\end{proposition}
Proof. We give the details for the maps $s \mapsto \mathcal T_s^{*-}(x, y)$, the proof goes along the same lines for $s \mapsto \mathcal T_s^{+}(x, y)$.
 Let $\Omega$ be the open neighborhood of $\overline{B(0, 1)}\setminus\{1\}$ given by Proposition \ref{localexpansion} and fix $\delta >0$ such that (\ref{mathfrakTx}), (\ref{local pour mathfrakT(x)}) and (\ref{uniformO}) hold. In particular,   we know that  for any $x, y \in \mathbb Z^-$, the function  $s \mapsto \mathcal T^{*-}_s(x, y)= \mathfrak T^{*-}(s\vert y-x)$ is  analytic on $\Omega$
and has   the local expansion, for $s \in {\mathcal O}_\delta(1)$ 
 $$
 \mathcal T^{*-}_s(x, y) =  {\mathcal T}^{*-}(x, y) + \sqrt{1-s}\ \widetilde{\mathcal T}^{*-}(x, y) +(1-s) \ {\bf O}_s(x, y)
 $$
whose coefficients are the ones  given in the statement of the proposition and $s \mapsto{\bf O}(x, y)$ is analytic in the variable $\sqrt{1-s}$; furthermore, the quantities
 $K^{\vert y-x\vert} \Big\vert \mathcal T^{+}_s(x, y)\Big\vert$ and $K^{\vert y-x\vert} \Big\vert {\bf O}_s(x, y)\Big\vert$ are bounded, uniformly  in  $x, y \in \mathbb Z^-$ and $s \in \Omega \cup {\mathcal O}_\delta(1)$.

\rightline{$\Box$}

\subsection{ The coefficient maps $s \mapsto \mathcal U_s^{*-}(x, y)$ and $s \mapsto \mathcal U_s^+(x, y)$ for $x, y \in \mathbb Z$}

 We   consider here the maps $s \mapsto  \mathcal U_s^{*-}(x, y)$ and $s \mapsto  \mathcal U_s^+(x, y)$.  Formally, the matrice $U_s^{*-} = \left(\mathcal U_s^{*-}(x, y)\right)_{x, y \in \mathbb Z}$ is the potential of $\mathcal T_s^{*-} = \left(\mathcal T_s^{*-}(x, y)\right)_{x, y \in \mathbb Z}$; since $\mathcal T_s^{*-} $  is strictly upper triangular, each $\mathcal U_s^{*-}(x, y)$ will be the combination by summations and products of finitely many coefficients $\mathcal T_s^{*-}(i, j), i, j \in \mathbb Z,$ and their regularity will thus be a direct consequence of the previous statement.  It will be a little more delicate for the coefficients of the matrice $\displaystyle \mathcal U_s^{+}= \sum_{n \geq 0} \left(\mathcal T_s^{+}\right)^n$ since the matrice $\mathcal T_s^{+}$ is upper triangular with non zero terms on the diagonal; we will mention  the adjustments we need  in this case.
One gets the
 
 \begin{proposition} \label{Uplusmoinss(x,y)}
There exists an open neighborhood $\Omega$ of $\overline{B(0, 1)}\setminus\{1\}$ such that, for any $x, y \in \mathbb Z^-$, the  functions  $s\mapsto \mathcal U^{*-}_s(x, y) $  
have an  analytic continuation to   $\Omega$; furthermore, for  $\delta>0$ small enough, these functions are analytic in the variable $\sqrt{1-s}$ 
 on   the set  $ {\mathcal O}_\delta(1)$    and   their  local expansions of order 1  are given by 
\begin{equation}\label{local pour Ustar-}
\mathcal U^{*-}_s(x, y) =  {\mathcal U}^{*-}(x, y) + \sqrt{1-s}\ \widetilde{\mathcal U}^{*-}(x, y) +(1-s) \ {\bf O}_s(x, y)
\end{equation}

where

$\bullet\quad $   
$\quad \displaystyle \mathcal U^{*-}(x, y)= U^{*-}(y-x)$

\vspace{2mm}

$\bullet\quad $   
$\quad \displaystyle \widetilde{\mathcal U}^{*-}(x, y)= -{\sqrt{2}\over \sigma}U^{*-}\Bigl(]y-x, 0]\Bigr) $

\vspace{2mm}

$\bullet\quad $   ${\bf O}_s(x, y)$ is analytic in the variable $\sqrt{1-s}$ and   bounded  for $s \in {\mathcal O}_\delta(1)$.

\vspace{2mm}

\noindent Similarly, for any $x, y \in \mathbb N_0$, the  functions  $s\mapsto \mathcal U^{+}_s(x, y) $  
have an  analytic continuation to   $\Omega$ and  these functions are analytic in the variable $\sqrt{1-s}$ 
 on   the set  $ {\mathcal O}_\delta(1)$    with the    local expansions of order 1  given by 
\begin{equation}\label{local pour Uplus}
\mathcal U^{+}_s(x, y) =  {\mathcal U}^{+}(x, y) + \sqrt{1-s}\ \widetilde{\mathcal U}^{+}(x, y) +(1-s) \ {\bf O}_s(x, y)
\end{equation}

where

$\bullet\quad \mathcal U^{+}(x, y)= U^{+}(y-x)$

\vspace{2mm}

$\bullet\quad \displaystyle \widetilde{\mathcal U}^{+}(x, y)= - {\sqrt{2}\over \sigma}U^{+}\Bigl([0, y-x]\Bigr) $

\vspace{2mm}

$\bullet\quad  {\bf O}_s(x, y)$ is analytic in the variable $\sqrt{1-s}$ and   bounded  for $s \in {\mathcal O}_\delta(1)$.

\vspace{2mm}

%

\end{proposition}
Proof.  Formally, one gets $\displaystyle \mathcal U^{*-}_s= \sum_{n \geq 0} 
\left( \mathcal T^{*-}_s\right)^n$; since the matrix  is strictly upper triangular, for any $x, y \in \mathbb Z^-$, one gets $\left( \mathcal T^{*-}_s\right)^n(x,y)=0$ for any $n>\vert x-y\vert$, so
\begin{equation}\label{sommefinie}
 \mathcal U^{*-}_s(x, y)=\sum_{n=0}^{\vert x-y\vert} \left( \mathcal T^{*-}_s\right)^n(x, y).
\end{equation}
The analyticity dependence,  for fixed $x, y \in \mathbb Z^-$, of the coefficients  $  \mathcal U^{*-}_s(x, y)$,  with respect to $s \in \Omega$ and $\sqrt{1-s}$ when $s \in {\mathcal O}_\delta(1)$,  immediately follows from  the previous Proposition.

Let us now establish the local expansion (\ref{local pour Ustar-});   for any fixed $x, y \in \mathbb Z^-$, one gets 
$$
 \mathcal U^{*-}_s(x, y)= \sum_{n=0}^{\vert x-y\vert}
 \Bigl(
 {\mathcal T}^{*-} + \sqrt{1-s}\ \widetilde{\mathcal T}^{*-} +(1-s) \ {\bf O}_s
 \Bigr)^n(x, y).
$$
The constant term  $\mathcal U^{*-}(x, y)$ is thus equal to $ \displaystyle 
\sum_{n=0}^{\vert x-y\vert}
 \Bigl(
 {\mathcal T}^{*-}
 \Bigr)^n(x, y)=\sum_{n=0}^{+\infty}
 \Bigl(
 {\mathcal T}^{*-}
 \Bigr)^n(x, y) ;$ on the other hand, the coefficient  corresponding to   $\sqrt{1-s}$ in this expansion is equal to 
 $$\widetilde{\mathcal U}^{*-}(x, y)= \sum_{n =0}^{\vert x-y\vert} \sum_{k=0}^{n-1}
 \Bigl(
 {\mathcal T}^{*-}  \Bigr)^k \widetilde{\mathcal T}^{*-}  \Bigl(
 {\mathcal T}^{*-}  \Bigr)^{n-k-1}(x, y).$$ 
 Inverting  the order of summations   and using the expression of $ \widetilde{\mathcal T}^{*-} $ in Proposition \ref{Ts(x,y)}, one gets
\begin{eqnarray*}
\widetilde{\mathcal U}^{*-}(x, y)&=&{\mathcal U}^{*-}\widetilde{\mathcal T}^{*-} {\mathcal U}^{*-}(x, y)\\
&=&
 -{{\sqrt 2}\over \sigma}\left(U^{*-}*\Bigl(\sum_{k \leq -1}\mu^{*-}\left( ]-\infty, k]\right)\ \delta_k \Bigr) * U^{*-} \right)(y-x)\\
 &=& -{{\sqrt 2}\over \sigma}
 U^{*-} \left(]y-x, 0]\right) \quad 
 \end{eqnarray*}
(to obtain the last equality, one  compute  the generating function  of the  measure
$$ U^{*-}*\Bigl(\sum_{k \leq -1}\mu^{*-}( ]-\infty, k])\ \delta_k \Bigr) * U^{*-},$$ it is equal to the one of the measure $U^{*-}*\lambda^{*-}$, and one concludes checking that 
 $$U^{*-}*\lambda^{*-}=\displaystyle 
\sum_{k \leq -1}\ U^{*-} \left(]k, 0]\right)\ \delta_k.)$$
 
 The proof goes along the same  lines for $\displaystyle  \mathcal U^{+}_s(x, y)= \sum_{n=0}^{+\infty}
 \Bigl(
 \mathcal T^{+}_s\Bigr)^n(x, y)$. Nevertheless, since $\mu^+(0)>0$, there are infinitely many terms in the sum;  for $s \in \Omega \cup {\mathcal O}_\delta(1)$,  one thus first sets $\mathcal T^+_s=\varepsilon_sI+T_s$
 $$
 $$
 with $\varepsilon_s:= \mathbb E\left[s^{\tau^+}1_{\{0\}}(S_{\tau^+})\right]$.
One gets $\delta_1 =\mu^+(0)\in ]0, 1[$, so $\vert \varepsilon_s\vert
< 1$ for    $\Omega$ and $\delta$ small enough. Since $I$ and $T_s$ commute and $T_s$ is strictly upper triangular, one may write, for any $x, y \in \mathbb N_0$ and $n \geq \vert x-y\vert, $
 \begin{eqnarray*}
 \Bigl(\mathcal T^+_s\Bigr)^n(x, y)
 &=& \sum_{k=0}^n \left(\begin{array}{c}n\\k\end{array}\right) \varepsilon_s^{n-k} T_s^k(x, y)\\
 &=& \sum_{k=0}^{\vert x-y\vert} \left(\begin{array}{c}n\\k\end{array}\right) \varepsilon_s^{n-k} T_s^k(x, y)
 \end{eqnarray*}
 so that
 \begin{eqnarray*}
\mathcal U^+_s(x, y)
 &=&  \sum_{n \geq 0}  \Bigl(\mathcal T^+_s\Bigr)^n(x, y)
 \\
 &=&\sum_{n =0}^{\vert x-y\vert}  \Bigl(\mathcal T^+_s\Bigr)^n(x, y)
 + \sum_{n >\vert x-y\vert} \sum_{k=0}^{\vert x-y\vert} \left(\begin{array}{c}n\\k\end{array}\right) \varepsilon_s^{n-k} T_s^k(x, y)\\
 &=&\sum_{n =0}^{\vert x-y\vert}  \Bigl(\mathcal T^+_s\Bigr)^n(x, y)
 + \sum_{k=0}^{\vert x-y\vert}  
 {1\over k!}\left(\sum_{n >\vert x-y\vert} n\cdots (n-k+1) \varepsilon_s^{n-k}\right) T_s^k(x, y)
 \end{eqnarray*} 
with  
 $s\mapsto \displaystyle \left(\sum_{n >\vert x-y\vert} n\cdots (n-k+1) \varepsilon_s^{n-k}\right) $ analytic on $\Omega$ and analytic in $\sqrt{1-s}$ on ${\mathcal O}_\delta(1)$. The analyticity of the map $s\mapsto \mathcal U_s^+(x, y)$ follows immediately; the computation of the coefficients of the local expansion
 (\ref{local pour Uplus}) goes along the same line than the ones of   (\ref{local pour Ustar-}).
 
%

\rightline{$\Box$}

\section{The centered reflected random walk}
Throughout this section, we will assume that hypotheses {\bf H} hold and that $\mu$ is centered. In this case, the radius of convergence of the generating functions $\mathfrak G(\cdot\vert x,  y), x , y \in \mathbb N_0,$ is equal to $1$. By Darboux's theorem, the asymptotic behavior of  the Taylor  coefficients  of these   generating functions is related to the type of their singularity near $s=1$; in the following subsection, we state some preparatory results.

 We denote by $\mathcal M$ the space  of infinite  matrices $M=(M(x,y))_{x, y \in \mathbb N_0}$; we will consider the elements of $\mathcal M$   as operators acting on the Banach space $(\mathbb C^{\mathbb N_0}, \vert\cdot\vert_\infty)$  and will thus  endowed $\mathcal M$ with the norm $\Vert \cdot \Vert_\infty $ defined by
$$\forall M \in \mathcal M \qquad \Vert M\Vert_\infty = \sup_{x \in \mathbb N_0} \sum_{y\in \mathbb N_0}\vert M(x, y)\vert.$$
Notice that this is the norm of $M$ considered as an operator acting on the Banach space $(\mathbb C^\mathbb N_0, \vert \cdot \vert_\infty)$ where $\vert \cdot \vert_\infty$ denotes the norm of the supremum.

As we have already seen, we will also  endow $\mathcal C^{\mathbb N_0}$ with the norm $\vert\cdot\vert_K$, with $K \in \mathcal K(1+\eta)$ for some $\eta>0$; the corresponding operator norm $\Vert\cdot\Vert_K$ on $\mathcal M$ will thus be defined by
$$\forall M \in \mathcal M \qquad \Vert M\Vert_K = \sup_{x \in \mathbb N_0} \sum_{y\in \mathbb N_0}{K(y)\over K(x)}\vert M(x, y)\vert.$$

\subsection{ The $\mathcal M$-valued map $s \mapsto {\mathcal T}_s$ and its potential $\mathcal U_s$}

Recall that the matrix $\mathcal T_s$ is the lower triangular with coefficients $\mathcal T_s(x, y), x, y \in \mathbb N_0$, given by
$$
\mathcal T_s(x, y)= \mathfrak T_s(y-x)= \mathbb E\left[ s^{\tau^{*-}}1_{\{y-x\}}(S_{\tau^{*-}})\right].
$$
The following statement is thus a direct consequence of Proposition \ref{localexpansion}:

\begin{proposition}\label{themapTs}
There exists an open neighborhood $\Omega$ of $\overline{B(0, 1)}\setminus\{1\}$ such that  the $\mathcal M$-valued function  $s\mapsto \mathcal T_s $  
has an  analytic continuation to   $\Omega$; furthermore, for  $\delta>0$ small enough, this function  is  analytic in the variable $\sqrt{1-s}$ 
 on   the set  $ {\mathcal O}_\delta(1)$    and  its   local expansions of order 1  in $(\mathcal M,  \Vert\cdot\Vert_\infty)$ is  given by 
\begin{equation}\label{local pour T(x,y)}
\mathcal T_s =  {\mathcal T} + \sqrt{1-s}\ \widetilde{\mathcal T} +(1-s) \ {\bf O}_s
\end{equation}

 where 
 
 $\bullet\quad $   $\mathcal T= \Bigl(\mathcal T(x, y)\Bigr)_{x, y \in \mathbb N_0}$ with 
$\quad \displaystyle \mathcal T(x, y)=
\left\{
\begin{array}{ll}
\mu^{*-}(y-x)& {\rm if\ } 0\leq y\leq x-1\\
0 & {\rm if\ } y\geq x
\end{array}
\right.$,
\vspace{2mm}

$\bullet\quad $   $\widetilde {\mathcal T}= \Bigl(\widetilde {\mathcal T}(x, y)\Bigr)_{x, y \in \mathbb N_0}$
 with 
$\quad \displaystyle  \widetilde {\mathcal T}(x, y)=
\left\{
\begin{array}{ll}
-{\sqrt{2}\over \sigma}
\mu^{*-}\Bigl(]-\infty, y-x]\Bigr)& {\rm if\ } 0\leq y\leq x-1\\
0 & {\rm if\ } y\geq x
\end{array}
\right.
$,

\vspace{2mm}

$\bullet\quad $   ${\bf O}_s$ is analytic in the variable $\sqrt{1-s}$ and  uniformly bounded in  $(\mathcal M, \Vert \cdot\Vert_\infty)$ for $s \in {\mathcal O}_\delta(1)$.

\end{proposition}
Proof. The regularity of each coefficient map $s \mapsto \mathcal T_s(x, y)$  may be proved as in Proposition  \ref{Ts(x,y)}; we thus focuse our attention on the analyticity of the $\mathcal M$-valued map $s \mapsto \mathcal T_s$.
By a classical result in the theory of vector valued analytic functions of the complex variable (see for instance \cite{C},  Theorem 9.13), it suffices to check that this property is true for the functions $s\mapsto  \mathcal T_s({\bf a})$ for any bounded sequence ${\bf a} = (a_i)_{i \geq 0} \in \mathbb C^\mathbb N_0$; to check this, we will use the fact that
any uniform limit  on some open set of analytic functions is analytic on this set.

Fix $N\geq 1$ and let $\mathcal T_{s,N}$ be the ``truncated'' matrix defined by
$$
\mathcal T_{s,N}(x, y)= 
\left\{
\begin{array}{ll}
\mathcal T_{s}(x, y) & {\rm if\ }  \max (x-N, 0)\leq y\leq x-1\\
0 & {\rm otherwise. }
\end{array}
\right.
$$
One gets 
$\displaystyle
\mathcal T_{s,N} ({\bf a})= \sum_{1}^{N} \mathfrak T_s^{*-}(-k)  {\bf a}^{(k)} 
\qquad {\rm with} \quad {\bf a}^{(k)}:= \underbrace{0, \cdots, 0}_{ k \ times},  a_0, a_1, \cdots, $ which implies that the
   $\mathcal M$-valued map $s\to \mathcal T_{s,N}$ is   analytic on $\Omega$  and analytic in the variable $\sqrt{1-s}$ on ${\mathcal O}_\delta(1)$. The same property holds for the map $s\to \mathcal T_s$ since, by (\ref{uniformO}), one gets
$$ \Vert\mathcal T_s-\mathcal T_{s,N}\Vert_\infty =\sup_{x \in \mathbb N_0}  \sum_{\vert y-x\vert >N}\vert \mathcal T_s(x, y)\vert \leq \sum_{\vert y-x\vert >N}{{\bf O}\over K^{\vert x-y\vert}}= {{\bf O} \over (K-1)K^N}\stackrel{N\to +\infty}{\longrightarrow} 0. $$

\rightline{$\Box$}

Let us now give sense to the matrix $(I-\mathcal T_s)^{-1}$; formally one may write  $$(I-\mathcal T_s)^{-1}=\mathcal U_s:= \sum_{k\geq 1} (\mathcal T_s)^k.$$
Since the matrices $\mathcal T_s$ are strictly lower triangular, one gets $\mathcal T_s^k(x, y)= 0$ for any $x, y \in \mathbb N_0$ and $k\geq \vert x-y\vert +1$; it follows that, for any $x, y \in \mathbb N_0$
\begin{equation} \label{defUs}
(I-\mathcal T_s)^{-1}(x, y)=\mathcal U_s(x, y) = \sum_{k=0}^{\vert x-y\vert}(\mathcal T_s)^k(x, y).
\end{equation}
The analyticity in the variable $s$ (resp. $\sqrt{1-s})$ on $\Omega$ (resp. on  ${\mathcal O}_\delta(1)$)  of each coefficient  $\mathcal U_s(x, y)$ follows by the previous fact and one may compute its local expansion near $s=1$. Nevertheless, this property does not hold in the Banach space $(\mathcal M,  \Vert\cdot\Vert_\infty)$, as can  be seen easily in the following statement (clearly, the matrices $\mathcal U $ and
$\widetilde{\mathcal U} $  which appear in (\ref{local pourUs(x,y)})  do not  belong  to the Banach space $(\mathcal M,  \Vert\cdot\Vert_\infty)$), we have in fact to consider a norm of the type $\Vert\cdot\Vert_K$ on $\mathbb C^{\mathbb N_0}$ to obtain a similar statement.  We may  state the following

\begin{proposition}\label{Us}
Fix $\eta>0$ and  $K\in {\mathcal K}(1+\eta)$. There exists an open neighborhood $\Omega$ of $\overline{B(0, 1)}\setminus\{1\}$ such that  the   function  $s\mapsto \mathcal U_s$  
has an  analytic continuation to   $\Omega$, with values in $(\mathcal M, \Vert\cdot\Vert)_K)$; furthermore, for  $\delta>0$ small enough, this function is  analytic in the variable $\sqrt{1-s}$ 
 on   the set  $ {\mathcal O}_\delta(1)$    and its   local expansion of order 1 in $(\mathcal M, \Vert\cdot\Vert_K)$ is  given by 
\begin{equation}\label{local pourUs(x,y)}
\mathcal U_s =  {\mathcal U}+ \sqrt{1-s}\ \widetilde{\mathcal U} +(1-s) \ {\bf O}_s
\end{equation}

 where

$\bullet\quad \mathcal U=\left(  \mathcal U(x, y) \right)_{x, y \in \mathbb N_0}$ with $\quad  \mathcal U(x, y) = 
\left\{
\begin{array}{ll}
U^{*-}(y-x)& {\rm if\ } 0\leq y\leq x\\
0 & {\rm if\ } y >x,
\end{array}
\right.
$,

\vspace{2mm}

$\bullet\quad  \widetilde{\mathcal U}=\left(   \widetilde{\mathcal U}(x, y) \right)_{x, y \in \mathbb N_0}$ with $\quad   \widetilde{\mathcal U}(x, y) = 
\left\{
\begin{array}{ll}
-{\sqrt{2}\over \sigma}
U^{*-}\Bigl(]y-x, 0]\Bigr)& {\rm if\ } 0\leq y\leq x-1\\
0 & {\rm if\ } y\geq x
\end{array}
\right.
$,
\vspace{2mm}

$\bullet\quad {\bf O}_s=\left( {\bf O}_s(x, y) \right)_{x, y \in \mathbb N_0}$
 is analytic in the variable $\sqrt{1-s}$  for $s \in {\mathcal O}_\delta(1)$   and uniformly bounded  in $(\mathcal M, \Vert \cdot\Vert_K)$.
\end{proposition}

\noindent Proof. Since $\Vert \mathcal T\Vert_\infty=1$, one may choose $\delta>0$ in such a way $\Vert \mathcal T_s\Vert_\infty\leq 1+{\eta\over 2}$ for any $s \in {\mathcal O}_\delta(1)$; it thus follows that, for such $s$,  any $x \in \mathbb N_0$ and $y\in \{0, \cdots, x-1\}$
\begin{equation}\label{majorUs}
\vert \mathcal U_s(x, y)\vert\leq \sum_{n=0}^{\vert x-y\vert}\Vert \mathcal T_s\Vert_\infty^k\leq (1+2/\eta) \left(1+{\eta\over 2}\right)^{\vert x-y\vert}.
\end{equation}
So, $\Vert \mathcal U_s \Vert_K<+\infty$ when $s \in {\mathcal O}_\delta(1)$ and  $K\in {\mathcal K}(1+\eta)$. To prove the analyticity of the function $s \mapsto \mathcal U_s$, we consider as above the truncated matrix $ \mathcal U_{s, N}$ and check, first that  for any ${\bf a}\in \mathbb C^{\mathbb N_0}$ the maps $s \mapsto  \mathcal U_{s, N}({\bf a})$ are analytic on $\Omega$ and analytic in the variable $\sqrt{1-s}$ on ${\mathcal O}_\delta(1)$,   and second that the sequence $( \mathcal U_{s, N})_{N\geq 1}$ converges to $ \mathcal U_{s}$ in $(\mathcal M, \Vert \cdot \Vert _K)$. The expansion (\ref{local pourUs(x,y)})
is a straightforward computation.

\rightline{$\Box$}

\noindent {\bf From now on, we fix a constant $\eta>0$ and a function $K$  in $\mathcal K(1+\eta)$.}

\subsection{The excursions $\mathcal E_s(\cdot, y)$ for $y \in \mathbb N_0$}

 The excursion  $\mathcal E_s$ before the first reflection  has been defined formally in (\ref{factor1bis})  as follows
$$\mathcal E_s= \Bigl(I-\mathcal  T_s\Bigr)^{-1}\mathcal U^+_s= \mathcal U_s\ \mathcal U^+_s.$$
The regularity  with respect to the parameter $s$  of the matrix coefficients  $ \mathcal U_s^+(x, y)$ and  the matrix $\mathcal U_s =\Bigl(I-\mathcal T_s\Bigr)^{-1}$  is well described  in Propositions \ref{Uplusmoinss(x,y)} and \ref{Us}. Each coefficient of $\mathcal E_s$ is a finite sum of products of coefficients of $\mathcal U_s$ and $\mathcal U^+_s$ so the regularity of  the map $s \mapsto  \mathcal E_s(x, y)$   will follow immediately;  the number of terms in this sum is equal to $\min(x, y)$, it thus increases with $x$ and $y$ and it is not easy  to obtain   some kind of uniformity with respect to these parameters. 
In fact, it will be sufficient to fix the arrival site $y$ and to describe the regularity of the $\mathbb C^\mathbb N_0$-valued map $s\mapsto (\mathcal E_s(x, y))_{x \in \mathbb N_0}$; to do this, we endow the space $\mathbb C^\mathbb N_0$  with the norm $\vert\cdot\vert_K$ defined in (\ref{normK}). 

\ 

We have the

\begin{proposition}\label{excursion}

There exists  an open neighborhood $\Omega$ of $\overline{B(0, 1)}\setminus\{1\}$ (depending on the function $K$) such that, for any $y \in \mathbb N_0$, the  functions    $s \mapsto \mathcal E_s(\cdot, y)$
have an  analytic continuation on $\Omega$    with values in the Banach space $(\mathbb C^\mathbb N_0, \vert \cdot\vert_K)$;  furthermore, for  $\delta>0$ small enough, these functions are  analytic in the variable $\sqrt{1-s}$ 
 on   the set  ${\mathcal O}_\delta(1)$ and  their local expansions of order 1 in $(\mathbb C^\mathbb N_0, \vert \cdot\vert _K)$ are  given by
\begin{equation}\label{local pour mathcalE}
\mathcal E_s(\cdot, y) = \mathcal E(\cdot, y)+ \sqrt{1-s}\   \widetilde{\mathcal E} (\cdot, y)+
 (1-s)\ {\bf O}_s(y)\end{equation}
  where 

$\bullet\quad $   $\mathcal E(\cdot, y)= \Bigl(I-\mathcal T\Bigr)^{-1}\mathcal U^+(\cdot, y)
  = \mathcal U\mathcal U^+(\cdot, y),$
  
  $\bullet\quad $   $\widetilde{\mathcal E}(\cdot, y)= \widetilde{ \mathcal U}\mathcal U^{+} (\cdot, y)+{\mathcal U}\widetilde{\mathcal U}^{+} (\cdot, y), $

$\bullet\quad $   ${\bf O}_s(y)$ is analytic in the variable $\sqrt{1-s}$ and  uniformly bounded in $(\mathbb C^\mathbb N_0, \vert \cdot\vert _K)$     for $s \in {\mathcal O}_\delta(1)$. 
\end{proposition}
Proof. Note that, for any $x \in \mathbb N_0$, one gets
$\displaystyle 
\mathcal E_s(x, y)= \sum_{z=0}^y\mathcal U_s(x, z)\mathcal U^+_s(z, y).
$
So, for $x$ fixed,  the conclusions above follows  from Propositions \ref{Uplusmoinss(x,y)} and \ref{coeffUs}; in particular, for any fixed $N\geq 1$, the  $\mathbb C^\mathbb N_0$-valued map $s\mapsto \Bigl(\mathcal E_{s, N}(x, y)\Bigr)$ defined by $\mathcal E_{s, N}(x, y)=\mathcal E_{s}(x, y)$ if $0\leq x\leq N$  and $\mathcal E_{s, N}(x, y)$ otherwise,  is analytic in $s\in \Omega$ and   $\sqrt{1-s}$ when $s \in {\mathcal O}_\delta(1)$, with values in the Banach space $(\mathbb C^\mathbb N_0, \vert \cdot\vert _K)$.  It is sufficient to check  that this sequence of vectors converges  to $\mathcal E_s(\cdot, y)_y$
in norm $\vert\cdot\vert_K$ for some suitable choice of $K>1$; by  (\ref{majorUs}), 
one gets  
$$
\Big\vert \mathcal E_s(x, y) \Big\vert \leq (y+1)(1+2/\eta) \left(1+{\eta\over 2}\right)^{x} \times \max_{0\leq z\leq y}\vert \mathcal U_s^+(z, y)\vert 
$$
so that 
$\displaystyle 
{\Big\vert \mathcal E_s(x, y) \Big\vert \over (1+\eta/2)^{x}}\leq {C_y\over (1+\eta/2)^x}$, for some constant $C_y>0$  depending only on $y$. Since $K\in \mathcal K(1+\eta)$,   one gets $\displaystyle \sup_{x\geq N} 
{\Big\vert \mathcal E_s(x, y) \Big\vert \over K(x)}\to 0$ as $N\to+\infty$; this proves that the sequence $\Bigl(\mathcal E_{s, N}(\cdot , y)\Bigr)_{N\geq 0} $  converges in $(\mathbb C^\mathbb N_0, \vert\cdot\vert_K)$ to $\mathcal E(\cdot, y)$ as $N\to+\infty$ and that $s \mapsto  \mathcal E_s(\cdot, y)$ is analytic. The local expansion (\ref{local pour mathcalE}) follows by a direct computation.

\rightline{$\Box$}

 \subsection{On the $\mathcal M$-valued map  $s \mapsto \mathcal R_s$}
 
 The matrices $\mathcal R_s$ which describe the dynamic of the space-time reflected process $({\bf r}_k, X_{{\bf r}_k})_{k \geq 0}$ is defined formally  in Section 2:
 $$
 \mathcal R_s=\Bigl(I-\mathcal T_s\Bigr)^{-1}\mathcal V_s= \mathcal U_s \mathcal V_s
 $$
 with  $\mathcal V_s=\Bigl(\mathcal V_s(x, y)\Bigr)_{x, y\in \mathbb N_0}$
where
$\displaystyle
\mathcal V_s(x, y):= \left\{
\begin{array}{ll}
0 &{\rm if} \ y=0\\
\mathfrak T^{*-}( s\vert -x-y ) &{\rm if} \ y\in \mathbb N^*
\end{array}
\right..
$
So, one first needs  to control the regularity of the map $s\mapsto \mathcal V_s$; as above, one gets the 
\begin{fact}\label{mathcalVs}
The function $s \mapsto \mathcal V_s$, with values in the Banach space $(\mathcal M, \Vert \cdot\Vert_K)$,  is analytic in $s$ on $\Omega$ and in the variable $\sqrt{1-s}$   on ${\mathcal O}_\delta(1)$; furthermore, it has    the following   local expansion  of order 1 near $s = 1$ \begin{equation} \label{localVs}
 \mathcal V_s = \mathcal V+ \sqrt{1-s}\   \widetilde{\mathcal V} +
 (1-s)\ {\bf O}_s
 \end{equation}
where 

\vspace{2mm}

$\bullet\quad $   $\mathcal V= \Bigl( \mathcal V(x, y)\Bigr)_{x,  y\in \mathbb N_0} $ with 
$\quad \displaystyle 
\mathcal V(x, y):= \left\{
\begin{array}{ll}
0 &{\rm if} \ y=0\\
\mu^{*-}(  -x-y ) &{\rm if} \ y\in \mathbb N^*.
\end{array}
\right.
$

\vspace{2mm}

$\bullet\quad $   $ \widetilde{\mathcal V}= \Bigl(  \widetilde{\mathcal V}(x, y)\Bigr)_{x,  y\in \mathbb N_0} $ with 
$\quad \displaystyle 
 \widetilde{\mathcal V}(x, y):= \left\{
\begin{array}{ll}
0 &{\rm if} \ y=0\\
- {\sqrt{2}\over \sigma} \mu^{*-}\Bigl(]-\infty, -x-y]\Bigr) &{\rm if} \ y\in \mathbb N^*.
\end{array}
\right.
$

\vspace{2mm}

$\bullet\quad $   ${\bf O}_s$ is analytic in the variable $\sqrt{1-s}$ and  uniformly bounded in  $(\mathcal M, \Vert\cdot\Vert_K)$ for $s \in {\mathcal O}_\delta(1)$.

\end{fact}
We now may  describe the regularity  of the map  $s \mapsto \mathcal R_s$: 
\begin{proposition}\label{analyticRs}
 The function  $s \mapsto \mathcal R_s$ has an  analytic continuation to an open neighborhood of $\overline{B(0, 1)}\setminus\{1\}$   with values in the Banach space $(\mathcal M, \Vert \cdot\Vert_K)$;  furthermore, for  $\delta>0$ small enough, this  function is analytic in the variable $\sqrt{1-s}$ 
 on   the set  ${\mathcal O}_\delta(1)$ and  its local expansion of order 1 in $(\mathcal M, \Vert \cdot\Vert_\infty)$ is  given by

 \begin{equation}\label{local pour mathcalR}
\mathcal R_s = \mathcal R + \sqrt{1-s}\    \widetilde{\mathcal R} +
 (1-s)\ {\bf O}_s\end{equation}
where

  $\bullet\quad $   $ \widetilde{\mathcal R}= \widetilde{\mathcal U}\mathcal V +\mathcal U\widetilde{\mathcal V}. $

$\bullet\quad $   ${\bf O}_s$ is analytic in the variable $\sqrt{1-s}$ and  uniformly bounded in $(\mathcal M, \Vert \cdot\Vert_\infty)$ for $s \in {\mathcal O}_\delta(1)$. 
\end{proposition}
Proof. The analyticity of this function with respect to the variables $s$ or $\sqrt{1-s}$ is clear by 
Proposition \ref{Us} and Fact \ref{mathcalVs} and one may write, for $s\in {\mathcal O}_\delta(1),$
\begin{eqnarray*}
\mathcal R_s &=& \Bigl(I-\mathcal  T_s\Bigr)^{-1} \mathcal V_s\\
&=& \mathcal U_s \mathcal V_s\\
&=&\Bigl( \mathcal U+\sqrt{1-s} \ \widetilde{\mathcal U}+(1-s) {\bf O}_s\Bigr)
\Bigl( \mathcal U+\sqrt{1-s} \ \widetilde{\mathcal U}+(1-s) {\bf O}_s\Bigr)\\
&=&  \mathcal U  \mathcal V+\sqrt{1-s}\Bigl(\ \widetilde{\mathcal U}\mathcal V+\mathcal U\widetilde{\mathcal V}\Bigr)+(1-s) {\bf O}_s. 
\end{eqnarray*}
 
\rightline{$\Box$}

A direct computation gives in particular
\begin{equation}\label{formulemathcalE}
 \mathcal E(x, y)= \sum_{k=0}^{\min(x, y)}U^{*-}(k-x)U^+
(y-k)\end{equation}
and
\begin{equation}\label{formulemathcalR'}
\widetilde{\mathcal R}(x, y)= \mathcal A(x, y)+\mathcal B(x, y)
\end{equation}
with
\begin{eqnarray*}  
\mathcal A(x, y)&:=& \left\{
\begin{array}{ll}
0 &{\rm if} \ x=0 \ {\rm or\ } y=0\\
\displaystyle  -{\sqrt{2}\over \sigma}\sum_{k=0}^{x-1}U^{*-}\Bigl(] k-x, 0]\Bigr)\mu^{*-}(-k-y)& {\rm otherwise}, 
\end{array}
\right.\\
\end{eqnarray*} 
and
\begin{eqnarray*}  
\mathcal B(x, y)&:=& \left\{
\begin{array}{ll}
0 &{\rm if} \  y=0\\
\displaystyle  -{\sqrt{2}\over \sigma}\sum_{k=0}^{x}U^{*-}(k-x)\mu^{*-}\Bigl(]-\infty, -k-y]\Bigr)& {\rm otherwise}. 
\end{array}
\right.
 \end{eqnarray*}

\subsection{ On the spectrum of $\mathcal R_s$ and its resolvant $\Bigl(I-{\mathcal R}_s\Bigr)^{-1}$} 

The question is more delicate  in the centered case since the spectral radius of $\mathcal R$ is equal to $1$ (we will see in the next Section that it is $<1$ in the non centered case, which simplify this step).

\subsubsection{The spectrum of $\mathcal R_s$ for $\vert s\vert = 1$ and $s \neq 1$}
{Using  Property \ref{quasicompact} , we first  control the spectral radius of the $\mathcal R_s$ for $s \neq 1$; indeed, we may control the norm of $\mathcal R_s^2$:
\begin{fact} \label{ spectralradius}
  For $\vert s\vert =1$ and $s\neq 1$ one gets $\Vert \mathcal R_s^2\Vert_K <1$; in particular, the spectral radius of $\mathcal R_s$ on $(\mathbb C^\mathbb N_0, \Vert \cdot \Vert_K)$ is $<1$.
\end{fact}
Proof. 
Fix $s \in \mathbb C\setminus\{1\}$  of modulus $1$;  by strict convexity, for any $w \in \mathbb N_0$ and $y \in \mathbb N^*$, there exists $\rho_{ w, y} \in ]0, 1[$, depending also on $s$,  such that $\vert \mathcal R_s(w, y)\vert \leq \rho_{ w, y}  \mathcal R(w, y)$;
on the other hand, by Property \ref{quasicompact},   we may choose $\epsilon >0$ and a finite set $F\subset \mathbb N_0$ such that, for any $x \in \mathbb N_0,$  
$$\mathcal R (x, F):= \sum_{w\in F}{\mathcal R}(x, w) \geq \epsilon.$$
 For any $y\in \mathbb N_0$, we set $\displaystyle \rho_y:= \max_{w\in F} \rho_{w, y}$; since $F$ is finite, one gets $\rho_y\in ]0, 1[$. 
 
Consequently, for any $x \in \mathbb N_0$ 
$$
 \Big\vert  \mathcal R_s^2K(x)\Big\vert 
\leq 
\sum_{w \in \mathbb N^*} \sum_{y \in \mathbb N^*}
 \mathcal R (x, w) \times \Big\vert \mathcal R_s(w, y)\Big\vert  K(y)
\leq \mathcal S_1(s\vert x)+ \mathcal S_2 (s\vert x)$$
with 
\begin{eqnarray*}
\mathcal S_1(s\vert x)&:= & \sum_{w \in F} \sum_{y \in \mathbb N^*}
  \mathcal R(x, w)  \times \Big\vert \mathcal R_s(w, y)\Big\vert  K(y)\\
  \mathcal S_2(s\vert x) &:= &
  \sum_{w \notin F} \sum_{y \in \mathbb N^*}
 \mathcal R(x, w) \times \Big\vert \mathcal R_s(w, y)\Big\vert   K(y).
 \end{eqnarray*}

One gets 
$$
\mathcal S_1(s\vert x)\leq  \sum_{w \in F} \mathcal R(x, w)  \sum_{y \in \mathbb N^*}
 \rho_y \mathcal R(w, y)K(y)\leq\rho  \mathcal R(x, F) $$
with $\displaystyle \rho := \max_{w \in F}\sum_{y \in \mathbb N^*}\rho_y \mathcal R(w, y)K(y) \in ]0, 1[.$ 

On the other hand
$
\displaystyle \mathcal S_2(s\vert x)\leq    \mathcal R(x, \mathbb N^*\setminus F)= 1-\mathcal R(x, F).$
Finally, since $K\geq1$, one gets 
$$
 {\Big\vert  \mathcal R_s^2K(x)\Big\vert \over K(x)}
  \leq  \Bigl( \rho  \mathcal R(x, F)+ 1-\mathcal R(x, F)\Bigr)\leq 1-(1-\rho)\epsilon <1,
$$
and the lemma follows.

\rightline {$\Box$}

Since the map $s \mapsto \mathcal R_s$ is analytic on the set $\{s \in \mathbb C/\vert s\vert<1+\delta\}\setminus[1, 1+\delta[$,  the same property holds for the map $s \mapsto (I-\mathcal R_s)^{-1}$ on a neighborhood of $\{s \in \mathbb C/\vert s\vert\leq 1 \}\setminus\{1\}.$

\subsubsection{Perturbation theory and spectrum of $\mathcal R_s$ for $s$ closed to $1$}
We now focus our attention on $s$ closed to   $1$. By Property \ref{quasicompact},  we know that the operator $\mathcal R$ may be decomposed as follows on $(\mathbb C^{\mathbb N_0}, \Vert\cdot\Vert_K)$
$$
\mathcal R= \pi+\mathcal Q
$$
where
 
$\bullet \  \pi$ is the rank one projector, on the space $\mathbb C\cdot {\bf h}$ generated by the sequence  ${\bf h}\in \mathbb C^{\mathbb N_0}$ whose terms are all equal to $1$, defined by
$$
{\bf a}=(a_k)_{k \geq 0}  \mapsto \Bigl(\sum_{i\geq1}\nu_{\bf r}(k)a_k \Bigr) {\bf h}\quad ^(\footnote{notice that   $\nu_{\bf r}({\bf h})=1$ since $\nu_{\bf r}$ is a probability measure on $\mathbb N_0$ }^),$$
 
 $\bullet \  \mathcal Q$ is a bounded operator on $(\mathbb C^{\mathbb N_0}, \Vert\cdot\Vert_K)$ with spectral radius $<1$,
 
 $\bullet\ \pi\circ \mathcal Q= \mathcal Q\circ \pi=0$.
 
Recall that  the map $s \mapsto \mathcal R_s$ is continuous on ${\mathcal O}_\delta(1)$ and, more precisely, that $s \mapsto \displaystyle { \mathcal R_s-\mathcal R \over \sqrt{1-s}}$ bounded on this set. By perturbation theory,  for   $s\in {\mathcal O}_\delta(1)$ with $\delta$ small enough, the  operator $\mathcal R_s$ admits  a similar spectral decomposition as above ; namely, one gets
\begin{equation}\label{spectraldecomposition}
\forall s \in {\mathcal O}_\delta(1) \qquad \mathcal R_s= \lambda_s\pi_s +\mathcal Q_s
\end{equation}
with 

 $\bullet \ \lambda_s$ is the dominant eigenvalue of $\mathcal R_s$, with corresponding eigenvector ${\bf h}_s$, normalized in such a way that $\nu_{\bf r}({\bf h}_s)=1$,

$\bullet \ \pi_s$ is a rank one projector on the space $\mathbb C\cdot {\bf h}_s$, 
 
 $\bullet \  \mathcal Q_s$ is a bounded operator on $(\mathbb C^{\mathbb N_0}, \Vert\cdot\Vert_K)$ with spectral radius $\leq \rho_\delta$ for some $\rho_\delta <1$,
 
 $\bullet\ \pi_s\circ \mathcal Q_s= \mathcal Q_s\circ \pi_s=0$.

\noindent Furthermore, the maps $\displaystyle s\mapsto {\lambda_s-1\over \sqrt{1-s}}, s\mapsto {\pi_s-\pi\over \sqrt{1-s}}, s\mapsto {{\bf h}_s-{\bf h}\over \sqrt{1-s}}$ and $\displaystyle s\mapsto {\mathcal Q_s-\mathcal Q\over \sqrt{1-s}}$ are bounded on ${\mathcal O}_\delta(1)$.
We may in fact precise  the local behavior of the map $s\mapsto \lambda_s$;  by the above decomposition and Proposition \ref{analyticRs}, one gets, for $s \in {\mathcal O}_\delta(1)$,
\begin{eqnarray*}
\lambda_s&=&\nu_{\bf r}(\mathcal R_s {\bf h}) +\nu_{\bf r}\Bigl((\mathcal R_s-\mathcal R)({\bf h}_s-{\bf h})\Bigr)\\
&=&1+\sqrt{1-s} \ \ \nu_{\bf r}(\widetilde{\mathcal R} {\bf h})+  (1-s)O(s)
\end{eqnarray*}
with $O(s)$ bounded on ${\mathcal O}_\delta(1)$. Since $\nu_{\bf r}(\widetilde{\mathcal R} {\bf h})\neq 0$, the operator $I-\mathcal R_s$ is invertible when $s \in {\mathcal O}_\delta(1)$ and $\delta$ small enough, with inverse
$$
\Bigl(I-\mathcal R_s\Bigr)^{-1}={1\over 1-\lambda_s}\ \pi_s + \Bigl(I-\mathcal Q_s\Bigr)^{-1}.
$$
We have thus obtained the   following
\begin{fact}\label{perturb-resolvante} For $\delta>0$ small enough, the function $s \mapsto \Bigl(I-\mathcal R_s\Bigr)^{-1}$ admits on ${\mathcal O}_\delta(1)$ the following  local expansion of order $1$ with values in $(\mathcal M, \Vert\cdot\Vert_K)$ 
\begin{equation}
\Bigl(I-\mathcal R_s\Bigr)^{-1}=-{1\over \sqrt{1-s}\times \nu_{\bf r}\Bigl(\widetilde{\mathcal R} {\bf h}\Bigr)}\  \pi+{\bf O}_s
\end{equation}
where ${\bf O}_s$ is analytic in the variable $\sqrt{1-s}$ and uniformly bounded in $(\mathcal M, \Vert\cdot\Vert_K).$
\end{fact}

\subsection{The return probabilities in the centered case: proof of the main theorem}
We use here the identity 
$\displaystyle 
\mathcal G_s=\Bigl(I-\mathcal R_s\Bigr)^{-1}\mathcal E_s
$
given in the introduction.
By  Proposition  \ref{analyticRs} and Fact \ref{ spectralradius}, for any fixed $y \in \mathbb N_0$, 
the function $s \mapsto \mathcal G_s(\cdot, y)$ is analytic on a neigborhood  of $\overline{B(0, 1)}\setminus\{1\}$. Furthermore, for $\delta>0$ small enough and $s \in {\mathcal O}_\delta(1)$, one may write, using (\ref{local pour mathcalE}) and (\ref{local pour mathcalR})
$$
\mathcal G_s(\cdot, y)= -{\nu_{\bf r}(\mathcal E(\cdot, y))\over \nu_{\bf r}\Bigl(\widetilde{\mathcal R} {\bf h}\Bigr) }\times {1\over \sqrt{1-s}}+{\bf O}_s
$$
with
$\nu_{\bf r}, \mathcal E(\cdot, y)$ and $\widetilde{\mathcal R}$ given respectively by formulae (\ref{stationary}), (\ref{formulemathcalE}) and (\ref{formulemathcalR'}) and $s \mapsto {\bf O}_s$  analytic on ${\mathcal O}_\delta(1)$ in the variable $\sqrt{1-s}$ and uniformly bounded in $(\mathcal M, \Vert \cdot\Vert_K)$.

We may thus apply Darboux's theorem \ref{darboux} with $R=1, \alpha =-{1\over 2}$ (and so   $\Gamma(-\alpha)=\sqrt{\pi})$)  and $\displaystyle \mathfrak{A}(1)= -{\nu_{\bf r}(\mathcal E(\cdot, y))\over \nu_{\bf r}\Bigl(\widetilde{\mathcal R} {\bf h}\Bigr) }>0$. One gets, for all $x, y \in \mathbb N_0$
\begin{equation}\label{constantCy}
\mathbb P_x[X_n=y]\sim{C_y\over \sqrt{n}}\quad \mbox{\rm with} \quad  C_y= - {1\over \sqrt{\pi}}\times{\nu_{\bf r}(\mathcal E(\cdot, y))\over \nu_{\bf r}\Bigl(\widetilde{\mathcal R} {\bf h}\Bigr) }>0.
\end{equation}
\rightline{$\Box$}

\section{The non centered random walk}

We  assume here $\mathbb E[Y_i]>0$ and use a standard argument in probability theory, called sometimes ''relativisation procedure'', to reduce the question  to the centered case.

\subsection{The relativisation principle and its consequences}
For any $r>0$, we denote by $\mu_r$ the probability measure defined  on $\mathbb Z$ by
$$\forall n \in \mathbb Z \qquad  \mu_r(n)= {1\over \hat{\mu}(r)} r^n \mu(n).$$

 Note that  for any $k \geq 0$ one gets $(\mu^{*k})_r=(\mu_r)^{*k}$ and that the generating function $\hat{\mu}_r$ is related to the one  of $\mu$ by the following identity
 $\displaystyle 
 \forall z \in \mathbb C \quad \hat{\mu}_r(z):= { \hat{\mu}(rz)\over \hat{\mu}(r)}.
 $

Notice that waiting times $\tau^{*-}$ and $\tau^+$ are defined on the  space $(\Omega, \mathcal T)$, with values in $\mathbb N_0\cup\{+\infty\}$, independently on the measure $\mu_r$ we choose; they are both a.s. finite if and only if  $\mu_r$ is centered, i.e. $r=r_0$ (see Section 3.1 for the notations). 

Throughout  this section, we will denote $\mathbb P^{\circ}$ the probability on $(\Omega, \mathcal T)$ which ensures that the $Y_i$ are i.i.d. with law $\mu_{r_0}$; the expectation with respect to $\mathbb P^{\circ}$ is denoted $\mathbb E^{\circ}$. We set $\rho_\circ= \hat{\mu}(r_0)$ and $R_\circ=1/\rho_\circ$; one gets  $R_\circ \in ]1, +\infty[$. Notice that the variables $Y_i$ have common law $\mu_{r_0}$  under $\mathbb P^\circ$, they  are in particular centered; we may thus apply the results of the previous section  when we refer to this probability measure on $(\Omega, \mathcal T)$.

We have the classical following 

\begin{fact}
Let $n\geq 1$ and $\Phi: \mathbb R^{n+1} \to \mathbb C$ a bounded Borel function; then, one gets
$$
\E\Bigl[\Phi(S_0, S_1, \cdots, S_n)\Bigr]= \rho_\circ^n\times  \E^\circ\Bigl[\Phi(S_0, S_1, \cdots, S_n)r_0^{-S_n}\Bigr].
$$
\end{fact}
As a direct consequence, for any $x, y \in \mathbb N_0$ and $s \in \mathbb C$, one gets, at least formally
$$
\mathcal E_s(x, y) = r_0^{x-y}\mathcal E^\circ_{\rho_\circ s}(x,y)\quad{\rm and}\quad 
\mathcal R_s(x, y) = r_0^{x+y}\mathcal R^\circ_{\rho_\circ s}(x,y).
$$
where we have set 
$\displaystyle
 \mathcal E^\circ_{ s}(x,y):= \sum_{n\geq 0} s^n\mathbb E^\circ_x[{\bf r}>n, X_n=y]
 $ 
 and 
 $\displaystyle
 \mathcal R^\circ_{ s}(x,y):= \sum_{n\geq 0} s^n\mathbb E^\circ_x[{\bf r}=n, X_n=y].
 $
 We may thus introduce the diagonal matrice $\Delta = (\Delta(x, y))_{x, y \in \mathbb N_0}$ defined by
 $\Delta(x, y)=0$ when $x \neq y$ and  $\Delta(x, x)=r_0^x$ for any $x\geq 0$; by the above, one gets formally
 $$
 \mathcal E_s=\Delta \mathcal E_s^\circ\Delta^{-1}
 \quad \mbox{\rm and} \quad 
  \mathcal R_s=\Delta \mathcal R_s^\circ\Delta.
 $$

In the sequel, we will add the exponent $\circ$ to the quantities $\mathcal U^+, \mathcal T, \mathcal U, \mathcal V$ defined in the previous section when $\mu$ was assume to be centered and considered here as variables defined on $(\Omega, \mathcal F, \mathbb P^\circ)$; with these notations, we will  have   $\mathcal E^\circ= \mathcal U^\circ \mathcal U^{\circ +}$, $\widetilde{\mathcal E}^\circ= \widetilde{\mathcal U}^\circ \mathcal U^{\circ +}+\mathcal U^\circ \widetilde{\mathcal U}^{\circ +}, \mathcal R^\circ= \mathcal U^\circ \mathcal V^{\circ }$ and
$\widetilde{\mathcal R}^\circ= \widetilde{\mathcal U}^\circ \mathcal V^{\circ }+\mathcal U^\circ \widetilde{\mathcal V}^{\circ}$.

 Combining Propositions \ref{excursion} and \ref{analyticRs}, we may thus state the

 \begin{proposition} \label{summaryNONCENTERED}There exist  a function $K$,  an open neighborhood $\Omega$ of $\overline{B(0, R_\circ)}\setminus\{R_\circ\}$ and $\delta>0$ small enough  such that, for any $y \in \mathbb N_0$, the  functions    $s \mapsto\Bigl( \mathcal E_s(x, y)\Bigr)_{x \in \mathbb N_0}$
have an  analytic continuation on $\Omega$    with values in the Banach space $(\mathbb C^\mathbb N_0, \vert \cdot\vert_K)$ and  are  analytic in the variable $\sqrt{R_\circ-s}$ 
 on   the set  ${\mathcal O}_\delta(R_\circ)):= B(R_\circ, \delta)\setminus[R_\circ, R_\circ+\delta)$, with the following local expansion  of order 1 in $(\mathbb C^\mathbb N_0, \vert \cdot\vert _K)$:   
 
\begin{equation}\label{local pour mathcalE}
\mathcal E_s  =  \mathcal E + \sqrt{R_\circ-s}\   \widetilde{\mathcal E} +
 (R_\circ-s)\    {\bf O}_s
 \end{equation}
  where

$\bullet\quad  \mathcal E=\Delta \mathcal E^\circ\Delta^{-1}$
  
$\bullet\quad  \widetilde{\mathcal E}=\sqrt{\rho_0} \ \Delta  \widetilde{\mathcal E}^\circ\Delta^{-1}$

$\bullet\quad $   $    {\bf O}_s=\left({\bf O}_s(x, y)\right)_{x\in \mathbb N_0}$ is analytic in the variable $\sqrt{R_\circ-s}$ and  uniformly bounded in $(\mathbb C^\mathbb N_0, \vert \cdot\vert _K)$     for $s \in {\mathcal O}_\delta(R_\circ)$. 

Similarly, the function  $s \mapsto \mathcal R_s=\Bigl(\mathcal R_s(x, y)\Bigr)_{x, y \in \mathbb N_0}$ has an  analytic continuation to $\Omega$,  with values in the Banach space $(\mathcal M, \Vert \cdot\Vert_K)$,  and is analytic in the variable $\sqrt{R_\circ-s}$ 
 on   the set  ${\mathcal O}_\delta(R_\circ)$ with the following local expansion of order 1 in $(\mathcal M, \Vert \cdot\Vert_K)$:  
 
 \begin{equation}\label{local pour mathcalR}
\mathcal R_s  =  \mathcal R + \sqrt{R_\circ-s}\   \widetilde{\mathcal R} +
 (R_\circ-s)\    {\bf O}_s
 \end{equation}
  where

$\bullet\quad  \mathcal R=\Delta \mathcal R^\circ\Delta^{}$
  
$\bullet\quad  \widetilde{\mathcal R}=\sqrt{\rho_0} \ \Delta  \widetilde{\mathcal R}^\circ\Delta^{}$

$\bullet\quad $   $    {\bf O}_s=\left({\bf O}_s(x, y)\right)_{x\in \mathbb N_0}$ is analytic in the variable $\sqrt{R_\circ-s}$ and  uniformly bounded in $(\mathbb C^\mathbb N_0, \vert \cdot\vert _K)$     for $s \in {\mathcal O}_\delta(R_\circ)$.

 \end{proposition}

 To prove the  main theorem in the non centered case, we will thus apply the same strategy than in the previous section. The proof simplifies in this case  since the operator $I-\mathcal R_s$ becomes invertible; namely, one gets the

\begin{fact} \label{resolvantNONCENTERED}For $K$ suitably choosen, $\delta>0$ small enough and any $s \in {\mathcal O}_\delta(R_\circ)$, the spectral radius of the  operator $\mathcal R_s$ on $(\mathcal M, \Vert\cdot\Vert_K)$ is $<1$.
\end{fact}
Proof.   It will be a direct consequence of the continuity  of the map $s \mapsto \mathcal R_s$ on $K_{\delta}(R_0)$ and the  inequality $\Vert \mathcal R_{R_\circ}\Vert_K <1$. Indeed, one gets, using the definition of $\mathcal R$ and setting $\phi:=1_{[-x, +\infty[}$
\begin{eqnarray*}
\Vert \mathcal R_{R_\circ}\Vert _K
&\leq&
\sup_{x \in \mathbb N_0} \sum_{y \geq 1}
\Bigl(\sum_{n \geq 0} 
R_\circ^n\mathbb P\left[\phi(S_1)\cdots \phi(S_{n-1})1_{\{-x-y\}}(S_n)\right]\Bigr)K(y)\\
&=& 
\sup_{x \in \mathbb N_0} \sum_{y \geq 1}
\Bigl(\sum_{n \geq 0} 
R_\circ^nr_0^{x+y}K(y)\mathbb P^\circ\left[\phi(S_1)\cdots \phi(S_{n-1})1_{\{-x-y\}}(S_n)\right]\Bigr)\\
&\leq & r_0
\sup_{x \in \mathbb N_0}  
\Bigl(\sum_{n \geq 0} 
R_\circ^n \mathbb P^\circ\left[ \phi(S_1)\cdots \phi(S_{n-1})(1-\phi)(S_n)\right]\Bigr)
\quad{\rm if}\quad r_0^{y-1}K(y)\leq 1 \quad {\rm for \ all} \quad y\geq 1\\
  &\leq& r_0
\end{eqnarray*}which achieves the proof, assuming that $y \mapsto r_0^{y-1}K(y)$ is $\leq 1$ on $\mathbb N_0$.

\rightline{$\Box$}

As a direct consequence, one may write
$$
\left(I-\mathcal R_s\right)^{-1} =\sum_{n\geq 0} \mathcal R_s^n.
$$
Furthermore, the map $s\mapsto \left(I-\mathcal R_s\right)^{-1}$ is analytic in the variable $s$ on $\Omega$ and analytic in the variable $\sqrt{R_\circ-s}$ on $O_\delta(R_\circ)$ and the local expansion near $R_\circ$ is
\begin{equation}\label{localexpRESOLVANTE}
\left(I-\mathcal R_s\right)^{-1}=\left(I-\mathcal R\right)^{-1}+
\sqrt{R_\circ-s}
\left(I-\mathcal R\right)^{-1}\widetilde{ \mathcal R}\left(I-\mathcal R\right)^{-1}+ (R_\circ-s)\cdots 
\end{equation}
 \subsection{The return probabilities in the non centered case: proof of the main theorem}We use here the identity 
$\displaystyle 
\mathcal G_s=\Bigl(I-\mathcal R_s\Bigr)^{-1}\mathcal E_s
$
given in the introduction.
By  Proposition  \ref{summaryNONCENTERED} and Fact \ref{resolvantNONCENTERED}, for any fixed $y \in \mathbb N_0$, 
the function $s \mapsto \mathcal G_s(\cdot, y)$ is analytic on a neigborhood  of $\overline{B(0, R_\circ)}\setminus\{R_\circ\}$. Furthermore, for $\delta>0$ small enough and $s \in {\mathcal O}_\delta(R_\circ)$, one may write, using Proposition  \ref{summaryNONCENTERED}  and the local expansion   (\ref{localexpRESOLVANTE})
$$
\mathcal G_s(\cdot, y)= \left(I-\mathcal R\right)^{-1}(\cdot, y) +
\sqrt{R_\circ-s}\Bigl(
\left(I-\mathcal R\right)^{-1}\widetilde{ \mathcal R}\left(I-\mathcal R\right)^{-1}\mathcal E(\cdot, y)+\left(I-\mathcal R\right)^{-1}\widetilde{ \mathcal E}(\cdot, y)\Bigr)+
 (R_\circ-s){\bf O}_s
$$
with
$s \mapsto {\bf O}_s$  analytic on ${\mathcal O}_\delta(R_\circ)$ in the variable $\sqrt{R_\circ-s}$ and uniformly bounded in $(\mathcal M, \Vert \cdot\Vert_K)$.

We may thus apply Darboux's theorem \ref{darboux} with $R=R_\circ, \alpha ={1\over 2}$ (and so   $\Gamma(-\alpha)=-2\sqrt{\pi})$)  and $\displaystyle \mathfrak{A}(R_\circ)=  
\left(I-\mathcal R\right)^{-1}\widetilde{ \mathcal R}\left(I-\mathcal R\right)^{-1}\mathcal E(\cdot, y)+\left(I-\mathcal R\right)^{-1}\widetilde{ \mathcal E}(\cdot, y)  <0$.

One gets, for all $x, y \in \mathbb N_0$
$$
\mathbb P_x[X_n=y]\sim C_{x,y}{\rho_\circ^n \over n^{3/2}}
$$
with $\rho_\circ= {1\over R_\circ}=\hat{\mu}(r_0)\in ]0, 1[$ and 
\begin{equation}\label{constantCx,y}
C_{x, y}= - {1\over 2\rho_\circ \sqrt{\pi}}\times\Bigl(\left(I-\mathcal R\right)^{-1}\widetilde{ \mathcal R}\left(I-\mathcal R\right)^{-1}\mathcal E(x, y)+\left(I-\mathcal R\right)^{-1}\widetilde{ \mathcal E}(x, y)\Bigr)>0
\end{equation}
where the matrices $\mathcal R, \widetilde{\mathcal R}, \mathcal E$ and $ \widetilde{\mathcal E}$ are given explicitely in  Proposition \ref{summaryNONCENTERED}.

\rightline{$\Box$}

\end{document}